\documentclass[a4paper,11pt]{amsart}
\usepackage{amsthm}
\usepackage{tikz}
\usepackage{epsfig}
\usepackage{amsfonts}
\usepackage{amsmath}
\usepackage{amssymb}
\usepackage{hyperref,latexsym}
\PassOptionsToPackage{unicode}{hyperref}

\usepackage[english]{babel}
\usepackage{amsmath}
\usepackage{amsfonts}
\usepackage{enumerate}
\usepackage{amssymb,mathtools,caption}
\usepackage{color}
\usepackage{mathrsfs}
\usepackage{xkeyval,xifthen,bm,twoopt}

\def\un{{\bf 1\!\!I}}

\newtheorem{theo}{Theorem}[section]
\newtheorem{lemm}[theo]{Lemma}
\newtheorem{coro}[theo]{Corollary}
\newtheorem{prop}[theo]{Proposition}
\newtheorem{rema}[theo]{Remark}
\newtheorem{defi}{Definition}

\theoremstyle{definition}

\newtheorem*{exampl*}{Examples}
\newtheorem*{conj*}{Conjecture}

\newtheorem*{hist*}{History for $p=2$}

\newcommand{\lbr}[1][{(}]{\left#1}
\newcommand{\rbr}[1][{)}]{\right#1}

 \everymath{\displaystyle}

\numberwithin{equation}{section}
\newcommand{\R}{\mathbb{R}}

\newcommand{\N}{\mathbb{N}}

\newcommand{\eqdef}{\stackrel{{\rm{def}}}{=}}

\newcommand{\ga} {\gamma}

\usepackage[a4paper]{geometry}
\newcommand{\de}{\delta}
\newcommand{\gpgt}[1][]{\ifthenelse{\isempty{#1}}{g \left( \frac{\ga}{2} \right)}{\frac{\ga}{#1}}}

\newcommand{\taylor}[4][]{
\ifthenelse{\equal{#1}{2}}{#2 \lbr #3 \rbr + #2'\lbr#3 \rbr \lbr #4\rbr + \frac{#2'' \lbr\zeta\rbr \lbr #4 \rbr^2}{2} }{}
\ifthenelse{\equal{#1}{3}}{#2 \lbr #3 \rbr + #2'\lbr#3 \rbr \lbr #4 \rbr + \frac{#2'' \lbr#3\rbr \lbr #4 \rbr^2}{2} + \frac{#2''' \lbr\zeta\rbr \lbr #4 \rbr^3}{3!} }{}
}

\newcommand{\taylorinv}[4][]{
\ifthenelse{\equal{#1}{2}}{#2\ \lbr #3 \rbr + \lbr #2\rbr '\lbr#3 \rbr\lbr #4\rbr + \frac{\lbr #2\rbr'' \lbr\zeta\rbr \lbr #4 \rbr^2}{2} }{}
\ifthenelse{\equal{#1}{3}}{#2 \lbr #3 \rbr + \lbr #2\rbr'\lbr#3 \rbr \lbr #4\rbr + \frac{\lbr #2\rbr'' \lbr#3\rbr \lbr #4 \rbr^2}{2} + \frac{\lbr #2\rbr''' \lbr\zeta\rbr \lbr #4 \rbr^3}{3!} }{}
}

\newcommand{\gpg}[2][]{
    \ifthenelse{\isempty{#1}}
    {\ifthenelse{\isempty{#2}}{g(\ga)}{g\lbr#2\rbr} }
    {\ifthenelse{\isempty{#2}}{g^{(#1)}(\ga)}{g^{(#1)}\lbr#2\rbr} }
    }

\newcommand{\gpu}[1][]{
		\ifthenelse{\isempty{#1}}%
		      {g(u)}
		      {g^{(#1)}(u)}
		      }

\newcommand{\gpgd}[2][]{
\ifthenelse{\isempty{#2}}{
		\ifthenelse{\isempty{#1}}{g(\ga)}{}
		\ifthenelse{\equal{#1}{1}}{g'(\ga)}{}
		\ifthenelse{\equal{#1}{2}}{g''(\ga)}{}
		\ifthenelse{\equal{#1}{3}}{g'''(\ga)}{}
		\ifthenelse{\equal{#1}{4}}{g''''(\ga)}{}}
		{
		\ifthenelse{\isempty{#1}}{g\lbr#2\rbr}{}
		\ifthenelse{\equal{#1}{1}}{g'\lbr#2\rbr}{}
		\ifthenelse{\equal{#1}{2}}{g''\lbr#2\rbr}{}
		\ifthenelse{\equal{#1}{3}}{g'''\lbr#2\rbr}{}
		\ifthenelse{\equal{#1}{4}}{g''''\lbr#2\rbr}{}}
		      }

\newcommand{\gpud}[1][]{
		\ifthenelse{\isempty{#1}}{g(u)}{}
		\ifthenelse{\equal{#1}{1}}{g'(u)}{}
		\ifthenelse{\equal{#1}{2}}{g''(u)}{}
		\ifthenelse{\equal{#1}{3}}{g'''(u)}{}
		}		      
		      
\newcommand{\vsig}{\varsigma}

\newcommand{\vu}[1][]{
		\ifthenelse{\isempty{#1}}%
		      {\vsig(u)}
		      {\vsig^{(#1)}(u)}
		      }

\newcommand{\vg}[1][]{
		\ifthenelse{\isempty{#1}}%
		      {\vsig(\ga)}
		      {\vsig^{(#1)}(\ga)}
		      }
		      
\newcommand{\gk}[1][]{
		\ifthenelse{\isempty{#1}}
		{\ga^q}
		{\ga^{q-#1}}
		}

\newcommand{\ginv}[2][]{
\ifthenelse{\isempty{#2}}{
		\ifthenelse{\isempty{#1}}{g^{-1} \lbr \ga \rbr}{}
		\ifthenelse{\equal{#1}{1}}{\lbr g^{-1}\rbr' \lbr \ga \rbr}{}
		\ifthenelse{\equal{#1}{2}}{\lbr g^{-1}\rbr'' \lbr \ga \rbr}{}
		\ifthenelse{\equal{#1}{3}}{\lbr g^{-1}\rbr''' \lbr \ga \rbr}{}
		}{
		\ifthenelse{\isempty{#1}}{g^{-1} \lbr #2 \rbr}{}
		\ifthenelse{\equal{#1}{1}}{\lbr g^{-1}\rbr' \lbr #2 \rbr}{}
		\ifthenelse{\equal{#1}{2}}{\lbr g^{-1}\rbr'' \lbr #2 \rbr}{}
		\ifthenelse{\equal{#1}{3}}{\lbr g^{-1}\rbr''' \lbr #2 \rbr}{}
		    }
		    }

\newcommand{\td}{T_{\de}}

\newcommand{\gtd}[1][]{
  \ifthenelse{\isempty{#1}}{y(\td)}{}
  \ifthenelse{\equal{#1}{1}}{g'(y(\td))}{}
  \ifthenelse{\equal{#1}{2}}{g''(y(\td))}{}
  \ifthenelse{\equal{#1}{3}}{g'''(y(\td))}{}
}

\newcommand{\ytd}[1][]{
  \ifthenelse{\isempty{#1}}{y(\td)}{}
  \ifthenelse{\equal{#1}{1}}{y'(\td)}{}
  \ifthenelse{\equal{#1}{2}}{y''(\td)}{}
  \ifthenelse{\equal{#1}{3}}{y'''(\td)}{}
}

\newcommand{\plap}[2][]{\ifthenelse{\isempty{#1}}{\lvert #2 \rvert^{n-2} #2}{\lbr \lvert #2 \rvert^{n-2} #2 \rbr}}

\newcommand{\gt}[2][]{\ifthenelse{\isempty{#2}}{
		\ifthenelse{\isempty{#1}}{g(y(t))}{}
		\ifthenelse{\equal{#1}{1}}{g'(y(t))}{}
		\ifthenelse{\equal{#1}{2}}{g''(y(t))}{}
		\ifthenelse{\equal{#1}{3}}{g'''(y(t))}{}
		}{
		\ifthenelse{\isempty{#1}}{g(y(#2))}{}
		\ifthenelse{\equal{#1}{1}}{g'(y(#2))}{}
		\ifthenelse{\equal{#1}{2}}{g''(y(#2))}{}
		\ifthenelse{\equal{#1}{3}}{g'''(y(#2))}{}
		    }
		    }

\newcommand{\psit}[2][]{\ifthenelse{\isempty{#2}}{
		\ifthenelse{\isempty{#1}}{\psi(y(t))}{}
		\ifthenelse{\equal{#1}{1}}{\psi'(y(t))}{}
		\ifthenelse{\equal{#1}{2}}{\psi''(y(t))}{}
		\ifthenelse{\equal{#1}{3}}{\psi'''(y(t))}{}
		}{
		\ifthenelse{\isempty{#1}}{\psi(y(#2))}{}
		\ifthenelse{\equal{#1}{1}}{\psi'(y(#2))}{}
		\ifthenelse{\equal{#1}{2}}{\psi''(y(#2))}{}
		\ifthenelse{\equal{#1}{3}}{\psi'''(y(#2))}{}
		    }
		    }

\newcommand{\nminn}[1][]{ \ifthenelse{\isempty{#1}}{ \lbr \frac{n-1}{n} \rbr}{\lbr \frac{n-1}{n} \frac{1}{#1}\rbr}} 

\newcommand{\funcpow}[4][y]{
\ifthenelse{\isempty{#2}}{\lbr #1 (#3) \rbr^{#4}}{}
\ifthenelse{\equal{#2}{1}}{\lbr #1' (#3) \rbr^{#4}}{}
\ifthenelse{\equal{#2}{2}}{\lbr #1'' (#3) \rbr^{#4}}{}
\ifthenelse{\equal{#2}{3}}{\lbr #1''' (#3) \rbr^{#4}}{}
\ifthenelse{\equal{#2}{4}}{\lbr #1'''' (#3) \rbr^{#4}}{}
\ifthenelse{\equal{#2}{5}}{\lbr #1''''' (#3) \rbr^{#4}}{}
\ifthenelse{\equal{#2}{6}}{\lbr #1'''''' (#3) \rbr^{#4}}{}
}

\usepackage{graphicx}
\usepackage{graphics}
\usepackage{epsfig}
\usepackage{amsmath}
\usepackage{amsfonts}
\usepackage{psfrag}

\begin{document}


\title[Fractional elliptic equation with singular nonlinearity]{Positive solutions to a fractional equation with singular nonlinearity}
\author{Adimurthi}
\address{T.I.F.R. CAM, P.B. No. 6503,\\
   Sharadanagar, Chikkabommasandra \\
   Bangalore 560065, India}
\email{adiadimurthi@gmail.com \and aditi@math.tifrbng.res.in}

\author{Jacques Giacomoni}
\address{Universit\'e de Pau et des Pays de l'Adour, LMAP (UMR CNRS 5142) Bat. IPRA,
   Avenue de l'Universit\'e \\
   F-64013 Pau, France}
\email{jacques.giacomoni@univ-pau.fr}

\author{Sanjiban Santra}
\address{Department of Basic Mathematics, Centro de Investigacione en Mathematicas,
   Guanajuato\\
  Mexico}
\email{sanjiban@cimat.mx}

%
%
%
%
%
%

%
%


\begin{abstract}
In this paper, we study the positive solutions to the following singular and non local elliptic problem posed in a bounded and smooth domain $\Omega\subset \R^N$, $N> 2s$: 
\begin{eqnarray*}
(P_\lambda)\left\{\begin{array}{lll}
&(-\Delta)^s u=\lambda(K(x)u^{-\delta}+f(u))\mbox{ in }\Omega\\
&u>0 \mbox{ in }\Omega\\
& u\equiv\, 0\mbox{ in }\R^N\backslash\Omega.
\end{array}\right.
\end{eqnarray*}
Here $0<s<1$, $\delta>0$, $\lambda>0$ and $f\,:\, \R^+\to\R^+$ is a positive $C^2$ function. $K\,:\, \Omega\to \R^+$ is a H\"older continuous function in $\Omega$ which behave as ${\rm dist}(x,\partial\Omega)^{-\beta}$ near the boundary with $0\leq \beta<2s$.

First, for any $\delta>0$ and for $\lambda>$ small enough, we prove the existence of solutions to $(P_\lambda)$. 
Next, for a suitable range of values of $\delta$, we show the existence of an unbounded connected branch of solutions to $(P_\lambda)$ emanating from the trivial solution at $\lambda=0$. For a certain class of nonlinearities $f$, we derive a global multiplicity result that extends results proved in \cite{peral-al}. To establish the results,  we prove new  properties which are of independent interest and deal with the behavior and H\"older regularity of solutions to $(P_\lambda)$.
\end{abstract}
\maketitle 
\section{Introduction}
Let $\Omega\subset \R^N$, $N> 2s$, be a bounded domain with boundary of class $C^{1,1}$. In this work, we study solutions to the Problem $(P_\lambda)$ above. Here $(-\Delta)^s$ is the fractional Laplace operator defined as 
\begin{equation*}
(-\Delta)^s u(x)=2 C(N,s){\rm P.V.}\int_{\R^N}\frac{u(x)-u(y)}{\vert x-y\vert^{N+2s}}\mathrm{d}y
\end{equation*}
where $\mathrm{P.V.}$ denotes the Cauchy principal value and $C(N,s)=\pi^{-\frac{N}{2}}2^{2s-1}s\frac{\Gamma(\frac{N+2s}{2})}{\Gamma(1-s)}$, $\Gamma$ being the Gamma function.

We assume that $0<s<1$, $\delta>0$, $\lambda\geq 0$, $K\in C_{{\rm loc}}^\nu(\Omega)$, $\nu\in (0,1)$, such that $\displaystyle\inf_{\Omega}K>0$ and satisfies for some $0\leq \beta<2s$ and $C_1, C_2>0$
\begin{equation}\label{eq1.1}
C_1d(x)^{-\beta}\leq K(x)\leq C_2 d(x)^{-\beta}, \quad \forall x\in \Omega
\end{equation}
where $d(x)\eqdef {\rm dist}(x,\partial\Omega)$.

Concerning $f$, we suppose the following conditions throughout the paper: 
\begin{itemize}
\item[{\bf (f1)}] $f\,:\,[0,\infty)\to\R$ is a positive $C^2$ function with $f(0)=0$;
\item[{\bf (f2)}] The function $g_x\,:\,t\to \frac{K(x)}{t^\delta}+f(t)$ is strictly convex on $(0,\infty)$  for any $x\in \Omega$;
\item[{\bf(f3)}] $\displaystyle\lim_{t\to\infty}\frac{f(t)}{t}=\infty$ and there exists $C>1$ such that $\displaystyle\liminf_{t\to\infty}\frac{f'(t)t}{f(t)}\geq C$.
\item[{\bf (f4)}] There exists $p\in \left(1, \frac{N+2s}{N-2s}\right)$ and $c>0$ such that $\displaystyle\lim_{t\to\infty}\frac{f(t)}{t^p}=c$.
\item[{\bf(f5)}] There exists $q\in \left(1, \frac{N+2s}{N-2s}\right)$ such that $\frac{tf'(t)}{f(t)}\leq q$, for any $t>0$.
\end{itemize}
The equation in $(P_\lambda)$ has intrinsic mathematical interest since in the local setting ($s=1$) it appears in several physical models like non newtonian flows in porous media, heterogeneous catalysts (see references \cite{DiMoOs}, \cite{FuMa}, \cite{GaJu}, \cite{GhRa} and the surveys \cite{GhRa2} and \cite{HeMa}). The fractional laplacian case has been investigated more recently in \cite{peral-al} and \cite{TuGiSe} where the existence and multiplicity of solutions have been proved by variational methods of mountain pass type and the non smooth analysis theory. In these two papers, the authors restrict to the case $f(u)=u^p$,  with $1<p\leq \frac{N+2s}{N-2s}$, $\beta=0$. Precisely, in \cite{peral-al}, the subcritical case (i.e. $1<p<\frac{N+2s}{N-2s}$) is considered. Existence of positive solutions are proved and a local multiplicity result is sketched for a certain range of $\delta$. In \cite{TuGiSe}, the critical case $p=\frac{N+2s}{N-2s}$ is dealt with and a global mul
 tiplicity result is proved for any $\delta>0$. The solutions have the form $u_\lambda=\underline{u}_\lambda+v_\lambda$ where $v_\lambda\in \tilde{H}^s(\Omega)$ and $\underline{u}_\lambda$ is the solution to the "pure singular" problem (see $(P_s)$ below),  i.e. $\underline{u}_\lambda$ satisfies: 
\begin{eqnarray*}
(P_s)\left\{\begin{array}{lll}
&(-\Delta)^s \underline{u}_\lambda=\lambda K(x){\underline{u}_\lambda}^{-\delta}\mbox{ in }\Omega,\\
&\underline{u}_\lambda>0 \mbox{ in }\Omega,\\
& \underline{u}_\lambda\equiv\, 0\mbox{ in }\R^N\backslash\Omega.
\end{array}\right.
\end{eqnarray*}

In the present paper, we investigate further Problem $(P_\lambda)$ for a larger class of nonlinearities $f$. We establish existence, multiplicity, asymptotic behaviour and regularity of solutions to $(P_\lambda)$. The multiplicity of solutions follows from the existence of global connected branch of solutions  that we prove by appealing the global bifurcation theory due to P.~H. Rabinowitz (see \cite{Ra-JFA}) in $\R^+\times C_{\phi_{\delta,\beta}}(\Omega)$, where $C_{\phi_{\delta,\beta}}(\Omega)$ is a closed subspace of $C_0(\overline{\Omega})$ weighted  by a suitable power of the distance to the boundary function.

In order to develop a bifurcation framework for Problem $(P_\lambda)$ in the $C_{\phi_{\delta,\beta}}(\Omega)$-setting, we need to prove new results about the behaviour of solutions to $(P_\lambda)$ and their H\"older-regularity in respect to the parameters $\delta$ and $\beta$. We point out that these results are new and also of independent interest. 

The asymptotic behaviour of solutions to $(P_\lambda)$ is also used to establish that the first eigenvalue of the linearized operator associated to the equation in $(P_\lambda)$, $\Lambda_1(\lambda)$, is well defined, principal and positive. Consequently from the implicit function theorem the branch of minimal solutions is smooth along the maximal interval $\lambda\in (0,\Lambda)$. 
 We also prove that $\Lambda_1(\lambda)$ is simple. From the Crandall-Rabinowitz local bifurcation result (see \cite{Cr-Ra}), we then deduce a local multiplicity result near $\lambda=\Lambda$.

Before stating precisely our main results, let us make some definitions. As in \cite{peral-al}, we adopt the following definition of (very) weak solutions:
\begin{defi}\label{weak_sol}
We say that $u\in L^1(\R^N)$, satisfying $u\equiv 0$ on $\R^N\backslash\Omega$, is a weak solution to $(P_\lambda)$ if $\inf_K u>0$  for any compact set $K\subset\Omega$ and for any $\phi\in \tau$,
\begin{equation}\label{form-var}
\int_{\Omega}u(-\Delta)^s\phi\,\mathrm{d}x=C(N,s)\int_{Q}\frac{(u(x)-u(y))(\phi(x)-\phi(y))}{\vert x-y\vert^{N+2s}}\,\mathrm{d}x\mathrm{d}y=\lambda\int_{\Omega}\left(\frac{K(x)}{u^\delta}+f(u)\right)\phi\,\mathrm{d}x
\end{equation}
where 
\begin{equation*}
\tau=\left\{\psi\,:\,\psi\,:\, \R^N\to\R, \mbox{ measurable and } (-\Delta)^s\psi\in L^\infty(\Omega),\; \psi\equiv 0\mbox{ on }\R^N\backslash\tilde{\Omega}, \tilde\Omega\subset\subset\Omega\right\}.
\end{equation*}
and

$Q=\R^{2N}\setminus(\mathcal C\Omega\times \mathcal C\Omega)$ and
 $\mathcal C\Omega := \R^N\setminus\Omega$. Note that $C^\infty_{\rm c}(\Omega)\subset\tau\subset L^\infty(\Omega)$.
\end{defi}
%
%
%
%
We then define the set of {\it classical} solutions to $(P_\lambda)$:
\begin{defi}\label{classical-sol}
Let 
\begin{equation*}
\mathcal{S}=\left\{(\lambda,u)\in \R^+\times C_0(\overline{\Omega})\,|\, u\mbox{ is a weak solution to } (P_\lambda)\right\}.
\end{equation*}
\end{defi}
\begin{rema}\label{about_class_sol}
Note that if $(\lambda, u)\in {\mathcal S}$, $(-\Delta)^su\in L^1_{\rm loc}(\Omega)$ and the equation in $(P_\lambda)$ is satisfied pointwise in $\Omega$. Furthermore, from the regularity theory for the fractional Laplacian, $u\in C^{1,\alpha}_{{\rm loc}}(\Omega)\cap C_0(\overline{\Omega})$. As we will see below, a comparison principle holds in the class of classical solutions.
\end{rema}
In the sequel, we will be interested to describe global and asymptotic properties of the set ${\mathcal S}$ in respect to the bifurcation parameter $\lambda$. In this regard, we make the following definition of an asymptotic bifurcation point:
\begin{defi}\label{def-asympt-bifur}
We call $\Lambda_a\in [0,\infty)$ an asymptotic bifurcation point for a subset $K$ of ${\mathcal S}$, if there exists a sequence $(\lambda_n,u_n)\in K$ such that $\lambda_n\to\Lambda_a$ and $\Vert u_n\Vert_{L^\infty(\Omega)}\to\infty$ as $n\to\infty$.
\end{defi}
Now, we define the space where solutions to $(P_\lambda)$ are setting: 
\begin{defi}
Given $\phi\in C_0(\overline{\Omega})$ such that $\phi>0$ in $\Omega$, define 
\begin{equation*}
C_\phi(\Omega)=\left\{u\in C_0(\Omega)\,|\,\exists c\geq 0\mbox{ such that }\vert u(x)\vert \leq c\phi(x), \;\forall x\in \Omega\right\}
\end{equation*}
\end{defi}
with the natural norm $\left\Vert \frac{u}{\phi}\right\Vert_{L^\infty(\Omega)}$ and the associated positive cone:
\begin{defi}
Define the following open  convex subset of $C_{\phi}(\Omega)$:
\begin{equation*}
C_{\phi}^+(\Omega)=\left\{u\in C_{\phi}(\Omega)\,|\, \displaystyle\inf_{x\in\Omega}\frac{u(x)}{\phi(x)}>0\right\}.
\end{equation*}
\end{defi}
Let $\phi_{1,s}$ be the first positive normalized eigenfunction ($\Vert\phi_{1,s}\Vert_{L^\infty(\Omega)}=1$) of $(-\Delta)^s$ in $\tilde{H}^s(\Omega)$ where

\begin{equation*}
\tilde{H}^s(\Omega)=\left\{u\in H^s(\R^N)\,|\, u= 0\mbox{ in }\mathcal C\Omega\right\}.
\end{equation*}
%
%
%
We recall that $\phi_{1,s}\in C^s(\R^N)$ and $\phi_{1,s}\in C_{d^s}^+(\Omega)$ (see for instance Proposition 1.1 and Theorem 1.2 in \cite{Ros-oton-serra-JMPA}). We then define  the function $\phi_{\delta,\beta}$ as follows:
\begin{eqnarray}\label{weighted-funct}
\phi_{\delta,\beta}=\displaystyle\left\{\begin{array}{lll}
& \phi_{1,s}\quad\mbox{ if } 0<\frac{\beta}{s}+\delta<1,\nonumber\\
& \phi_{1,s}\left(\ln\left(\frac{2}{\phi_{1,s}}\right)\right)^{\frac{1}{\delta+1}}\quad\mbox{ if } \frac{\beta}{s}+\delta=1,\nonumber\\
& \phi_{1,s}^{\frac{2s-\beta}{(\delta+1)s}}\quad\mbox{ if } \frac{\beta}{s}+\delta>1.
\end{array}\right.
\end{eqnarray}
We now give the statements of our all main results that we will prove in this paper. First, we deal with the pure singular problem that provides suitable subsolutions to $(P_\lambda)$:
\begin{eqnarray*}
(P_s)\displaystyle\left\{\begin{array}{lll}
& (-\Delta)^su=\frac{K(x)}{u^\delta}\mbox{ in }\Omega,\\
& u>0\mbox{ in }\Omega,\\
& u=0\mbox{ in }\R^N\backslash\Omega.
\end{array}\right.
\end{eqnarray*}
Considering $(P_s)$ we have the following results:
\begin{theo}\label{sing-prob}
\begin{itemize}
\item[i)] If $\frac{\beta}{s}+\delta\leq 1$, then there exists  a unique $u\in C_0(\overline{\Omega})$ classical solution to $(P_s)$. Furthermore, $u\in \tilde{H}^s(\Omega)\cap C_{\phi_{\delta,\beta}}^+(\Omega)$.
\item[ii)] If $\frac{\beta}{s}+ \delta>1$ with $\beta<2s$, then there exists $u\in C_{\phi_{\delta,\beta}}^+(\Omega)$ classical solution to $(P_s)$. Furthermore, $u\in \tilde{H}^s(\Omega)$ if and only if $2\beta+\delta(2s-1)<1+2s$ and in this case $u$ is the unique classical solution to $(P_s)$.
\item[iii)]If $s\geq 2s$, then there is no classical solution to $(P_s)$.
\end{itemize}
\end{theo}
\begin{rema}\label{dual_sense}
 If $2\beta+\delta(2s-1)<1+2s$, since $C^\infty_{\rm c}(\Omega)$ is dense in $\tilde{H}^s(\Omega)$ and by Hardy inequality (see \cite[Par. 3.2.6, Lem. 3.2.6.1,  p. 259]{Tr} or \cite[Corollary~1.4.4.10, p. 33]{Gr}), \eqref{form-var} (with $f=0$) is satisfied for any $\phi\in \tilde{H}^s(\Omega)$. Hence, any $\phi\in \tilde{H}^s(\Omega)$ can be used as a test function in  \eqref{form-var}.
\end{rema}
Concerning the H\"older regularity of weak solutions to $(P_s)$, we have the following:
\begin{theo}\label{sing-prob2}
\begin{itemize}
\item[i)] If $\frac{\beta}{s}+\delta< 1$, then the classical solution $u$ to $(P_s)$ belongs to $C^s(\R^N)$;
\item[ii)] If $\frac{\beta}{s}+\delta=1$, then the classical solution $u$ to $(P_s)$ belongs to $ C^{s-\epsilon}(\R^N)$ for any $\epsilon>0$ small enough;
\item[iii)]If $\frac{\beta}{s}+ \delta>1$ and $\beta<2s$, then any classical solution to $(P_s)$ in $C_{\phi_{\delta,\beta}}^+(\Omega)$ belongs to $C^{\frac{2s-\beta}{\delta+1}}(\R^N)$.
\end{itemize}
\end{theo}
\begin{rema}\label{regu-problemp}
Theorems~\ref{sing-prob} and \ref{sing-prob2} (see Proposition~\ref{minimal_branch}) hold for the following problem:
\begin{eqnarray*}
\left\{\begin{array}{lll}
&(-\Delta)^s u=\lambda(K(x)u^{-\delta}+g)\mbox{ in }\Omega\\
&u>0 \mbox{ in }\Omega\\
& u\equiv\, 0\mbox{ in }\R^N\backslash\Omega.
\end{array}\right.
\end{eqnarray*}
where $g\in L^\infty(\Omega)$ with similar proofs. Hence Theorems~\ref{sing-prob} and \ref{sing-prob2} are still valid for bounded solutions to $(P_\lambda)$.
\end{rema}
We now consider the problem $(P_\lambda)$. The next result shows the existence of a global branch of (classical) solutions to $(P_\lambda)$:
\begin{theo}\label{main-bifurcation}
Let $f$ satisfy conditions {\bf (f1)}-{\bf (f3)} and assume that $2\beta+\delta(2s-1)<1+2s$. Then, 
\begin{itemize}
\item[i)] There exists $\Lambda\in (0,+\infty)$ and $0<\gamma=\gamma(\beta,\delta)$ such that ${\mathcal S}\subset[0,\Lambda]\times \left(C^\gamma(\R^N)\cap \tilde{H}^s(\Omega)\cap C_{\phi_{\delta,\beta}}^+(\Omega)\right)$ with
\begin{eqnarray*}
\gamma=\displaystyle\left\{\begin{array}{lll}
& s\quad\mbox{ if }\,\frac{\beta}{s}+\delta<1,\\
& s-\epsilon \quad\mbox{ if }\,\frac{\beta}{s}+\delta=1, \;\forall\epsilon>0 \mbox{ small enough},\\
& \frac{2s-\beta}{\delta+1}\quad\mbox{ if }\,\frac{\beta}{s}+\delta>1.
\end{array}\right.
\end{eqnarray*}
\item[ii)] There exists a connected unbounded branch ${\mathcal C}$ of solutions to $(P_\lambda)$ in $\R^+\times C_0(\overline{\Omega})$, emanating from $(0,0)$ such that for any $\lambda\in (0,\Lambda)$, there exists $(\lambda, u_\lambda)\in {\mathcal C}$ with $u_\lambda$ being the minimal solution to $(P_\lambda)$. Furthermore, as $\lambda\to\Lambda^-$, $u_\lambda\to u_\Lambda$ in $\tilde{H}^s(\Omega)$, where $u_\Lambda$ is a weak solution to $(P_\Lambda)$.
\item[iii)] The curve $(0,\Lambda)\ni\lambda\mapsto u_\lambda\in C_0(\overline{\Omega})$ is of class $C^2$.
\item[iv)] ({\it Bending and local multiplicity near $\Lambda$})\\
If $u_\Lambda\in L^\infty$, then $\lambda=\Lambda$ is a bifurcation point, that is, there exists a unique $C^2$-curve $(\lambda(s),u(s))\in {\mathcal C}$, where the parameter $s$ varies in an open interval about the origin in $\R$, such that
\begin{equation}\label{bending-carac}
\lambda(0)=\Lambda,\; u(0)=u_\Lambda,\; \lambda'(0)=0,\; \lambda''(0)<0.
\end{equation}
\item[v)] ({\it Asymptotic bifurcation point})\\
${\mathcal C}$ admits an asymptotic bifurcation point $\Lambda_a$ satisfying $0\leq \Lambda_a\leq \Lambda$.
\end{itemize}
\end{theo}
Note that assertion (iv) in the above theorem implies that the connected branch bends to the left thereby creating at least two solutions in a left neighborhood of $\Lambda$. The next result provides a global multiplicity result for a class of functions $f$ including the case $f(u)=u^p$ with $1<p<\frac{N+2s}{N-2s}$.
\begin{theo}\label{bifurcation-example} Assume that $\beta=0$, $\delta(2s-1)<2s+1$ and that $f$ satisfies ${\bf (f1)}-{\bf (f5)}$. 
Then, for any $\lambda_0>0$,  ${\mathcal S}\cap \{\lambda\geq \lambda_0\}$ is uniformly bounded and the connected branch ${\mathcal C}$ of solutions to  $(P_\lambda)$ given in Theorem~\ref{main-bifurcation} admits one and only one bifurcation point $\Lambda_a=0$.
\end{theo}
The proof of Theorem~\ref{bifurcation-example} is based on uniform $L^\infty$-bound estimates. These estimates are established by blow-up technique together with the moving plane method. Furthermore, Theorem \ref{bifurcation-example} provides a global multiplicity result for the class of functions $f$ satisfying hypothesis {\bf(f1)}-{\bf (f5)}, i.e.
\begin{equation*}
\forall\lambda\in (0,\Lambda), \mbox{ there exist at least two distinct classical solutions to } \, (P_\lambda).
\end{equation*}
This result complements results in \cite{peral-al} and \cite{TuGiSe}.

In the following result, appealing the theory of analytic global bifurcation theory (see \cite{Buffoni-Toland-book}, \cite{Dancer}, \cite{GiPrWa}) and under additional  restrictions on $\delta$ and $\beta$, we prove the existence of a continuous and piecewise analytic curve of solutions to $(P_\lambda)$. We consider here the case where $f(u)=u^p$ with $1<p<\frac{N+2s}{N-2s}$ but it can be extended for more general $f$ with similar growth and as soon as the analyticity property for the operator $F$ defined below (see assertion (v) of Theorem~\ref{example-analytic}) holds.
\begin{theo}\label{example-analytic}
Let $f(t)=t^p$ for some $p\in \left(1, \frac{N+2s}{N-2s}\right)$. Assume that $\frac{\beta}{s}+\delta<1$. Then ${\mathcal C}$ contains  an unbounded set ${\mathcal A}$  which is globally parametrized by a continuous map:
\begin{equation*}
(0,\infty)\ni s\mapsto (\lambda(s),u(s))\in {\mathcal A}\subset {\mathcal S}.
\end{equation*}
Moreover, the following properties hold along the path ${\mathcal A}$:
\begin{itemize}
\item[i)] $(\lambda(s), u(s))\to (0,0)$ in $\R\times C_{\phi_{\delta,\beta}}(\Omega)$ as $s\to 0^+$.
\item[ii)] For some $s_0>0$, the portion of the path $\{(\lambda(s),u(s))\,:\, 0<s<s_0\}$ coincides with all the minimal solutions branch.
\item[iii)] $\Vert u(s)\Vert_{C_{\phi_{\delta,\beta}}}\to \infty$ as $s\to\infty$.
\item[iv)] ${\mathcal A}$ admits at least one asymptotic bifurcation point $\lambda_0\in[0,\Lambda]$, that is, there exists a sequence $(s_n)_{n\in\N}\subset (0,\infty)$ such that $\lambda(s_n)\to \lambda_0$ and $\Vert u(s_n)\Vert_{C_{\phi_{\delta,\beta}}}\to\infty$ as $n\to\infty$.
\item[v)] Let $F\, :\, \R^+\times C_{\phi_{\delta,\beta}}^+(\Omega)\to C_{\phi_{\delta,\beta}}^+(\Omega)$ defined by
\begin{equation*}
F(\lambda, u)= (-\Delta)^{-s}\left(\lambda(\frac{K(x)}{u^\delta}+u^p)\right)\quad\mbox{for any }(\lambda, u)\in \R^+\times C_{\phi_{\delta,\beta}}(\Omega)^+.
\end{equation*}
Then, $\left\{s\geq 0\,:\, \partial_u F(\lambda(s),u(s)) \mbox{ is not invertible }\right\}$ is a discrete set.
\item[vi)] ($\mathcal{A}$ is an ``analytic" path) At each of its points ${\mathcal A}$ has a local analytic re-parameterization in the following sense: For each $s^*\in (0,\infty)$ there exists a continuous, injective map $\rho^*\,:\, (-1,1)\to \R$ such that $\rho^*(0)=s^*$  and the re-parametrisation
\begin{eqnarray*}
 (-1,1) \ni t\to (\lambda(\rho^*(t)),u(\rho^*(t))) \in \mathcal{A} \mbox{ is analytic}.
\end{eqnarray*}
Furthermore, the map $s \mapsto \lambda(s)$ is injective in a right neighborhood of $s=0$ and for each $s^*>0$ there exists $\epsilon^*>0$ such that $\lambda$ is injective on $[s^*,s^*+\epsilon^*]$ and on $[s^*-\epsilon^*,s^*]$.
\item[vii)] ${\mathcal A}$ bends to the left of $\{\lambda = \Lambda\}$ at the point $(\Lambda, u_\Lambda)$.
 \end{itemize}
\end{theo}
\begin{rema}
If $K\equiv 1$, then from Theorem~\ref{bifurcation-example}, the asymptotic bifurcation point $\lambda_0=0$.
\end{rema}
In the proof of Theorem~\ref{main-bifurcation}, we will consider the compact nonlinear operator $A\,:\, \R^+\times C_0(\overline{\Omega})\mapsto C_0(\overline{\Omega})$ defined for any $(\lambda, v)\in \R^+\times C_0(\overline{\Omega})$ by:
\begin{equation*}
A(\lambda,v)=w\in \tilde{H}^s(\Omega)\mbox{ be the unique solution to }(-\Delta)^s w-\frac{\lambda K(x)}{w^\delta}=v, w\equiv 0\mbox{ on }\R^N\backslash\Omega.
\end{equation*}
We have clearly that
\begin{equation*}
(\lambda, u)\in {\mathcal S}\Leftrightarrow u=A(\lambda,\lambda f(u)).
\end{equation*}
The compactness of $A$ follows from the regularity result stated in assertion (i) of Theorem\ref{main-bifurcation} and is used to apply the global bifurcation theory of  Rabinowitz. 

Concerning the proof of Theorem~\ref{example-analytic}, we need the following result to get the compactness of the operator $F$ in $C_{\phi_{\delta,\beta}}(\Omega)$.
\begin{lemm}\label{Compactness-Lemma}
Let $(w_n)_{n\in\N}$ be a bounded sequence in $C_{\phi_{\delta,\beta}}(\Omega)$. Let $\Omega\ni x\to C(x)$ be a $C^{\gamma}_{\rm loc}$ positive function, with $\gamma\in (0,1)$,  such that $\displaystyle\sup_\Omega C(x)d(x)^{\beta+(\delta+1)s}<\infty$. Then,
\begin{equation*}
v_n=\left[(-\Delta)^s+C(x)I\right]^{-1} w_n\in C_{\phi_{\delta,\beta}}(\Omega)
\end{equation*}
and $(v_n)_{n\in\N}$ is relatively compact in $C_{\phi_{\delta,\beta}}(\Omega)$.
\end{lemm}
The analytic framework requires to prove additionally the analyticity of $F$ on $\R^+\times C_{\phi_{\delta,\beta}}(\Omega)$. This can be proved similarly as in \cite{Da2}.
\begin{rema}
From assertions (v) and (vi) in Theorem~\ref{example-analytic}, ${\mathcal A}$ has a nice structure made of analytic arcs parametrized by $\lambda$ and only a finite countable collection of singular points. 
\end{rema}
The outline of the paper is as follows. In Section~\ref{section1}, we prove Theorems~\ref{sing-prob} and \ref{sing-prob2}. In Section~\ref{section2}, we establish our main bifurcation result, that is Theorem~\ref{main-bifurcation}. Then, in Section~\ref{section3}, we deal with the special case of subcritical nonlinearities and prove Theorems~\ref{bifurcation-example} and \ref{example-analytic} together with Lemma~\ref{Compactness-Lemma}. Finally in the appendix, we prove the $C^2$-regularity of the operator $A$ involved in the proof of Theorem~\ref{main-bifurcation}.
\section{Pure singular problem $(P_s)$.}\label{section1}
In this section, we deal with the problem $(P_s)$ and prove Theorems~\ref{sing-prob} and \ref{sing-prob2}. We start with the proof of Theorem~\ref{sing-prob}.

\noindent{\bf Proof of Theorem~\ref{sing-prob}.}\\
We first prove the existence of classical solutions. We give two alternative proofs. Let us consider first the case $\frac{\beta}{s}+\delta<1$. In the spirit of the seminal work of Crandall, Rabinowitz, Tartar (see \cite{Cr-Ra-Ta}),  we introduce the following approximated problem:
\begin{eqnarray*}
(P_\epsilon)\left\{\begin{array}{lll}
& (-\Delta)^s u=\frac{K_\epsilon(x)}{(u+\epsilon)^\delta}\mbox{ in }\Omega,\\
&u>0\mbox{ in }\Omega,\\
& u\equiv 0\mbox{ in }\R^n\backslash\Omega
\end{array}\right.
\end{eqnarray*}
with $K_\epsilon(x)=\inf\left\{\frac{1}{\epsilon}, K(x)\right\}$. $(P_\epsilon)$ admits a unique solution  in $\tilde{H}^s(\Omega)$. Indeed, Let $\tilde{H}^{s}(\Omega)^+$ denote the positive cone of $\tilde{H}^s(\Omega)$. Let $E_\epsilon\,:\, \tilde{H}^{s}(\Omega)^+\to\R$ defined by:
\begin{equation*}
E_\epsilon(v)=\frac{1}{2}\int_{\R^N}((-\Delta)^{\frac{s}{2}} v)^2\,\mathrm{d}x-\int_{\Omega}\frac{K_\epsilon(x)(v+\epsilon)^{1-\delta}}{1-\delta}\,\mathrm{d}x
\end{equation*}
for any $v\in \tilde{H}^{s}(\Omega)^+$. It is easy to prove that $E_\epsilon$ is weakly lower semi-continuous, strictly convex and coercive on $\tilde{H}^{s}(\Omega)^+$. Furthermore, $\displaystyle\inf_{\tilde{H}^{s}(\Omega)^+}E_\epsilon<0$. Therefore, $E_\epsilon$ admits a unique global minimizer, $u_\epsilon\not\equiv\,0$, on $\tilde{H}^{s}(\Omega)^+$. Furthermore, we observe that for $c>0$ small enough,  $c\phi_{1,s}$ is a strict subsolution to $(P_\epsilon)$, independently of $\epsilon$. Indeed, for a constant $c>0$ small enough and independent of $\epsilon$, we have
\begin{equation*}
(-\Delta)^s(c\phi_{1,s})\frac{(c\phi_{1,s}+\epsilon)^\delta}{K_\epsilon(x)}=\lambda_{1,s}\frac{c\phi_{1,s}(c\phi_{1,s}+\epsilon)^\delta}{K_\epsilon(x)}\leq 1.
\end{equation*}
To derive the Euler-Lagrange equation, we prove that $u_\epsilon\geq c\phi_{1,s}$. For that, let us consider the convex function $\xi\,:\,[0,1]\to\R$ defined by
\begin{equation*}
\xi(t)=E_\epsilon(u_\epsilon+t(c\phi_{1,s}-u_\epsilon)^+).
\end{equation*}
By convexity of $\xi$ and since $u_\epsilon$ is a minimizer of $E_\epsilon$, $0\leq \xi'(0)\leq \xi'(1)$ and 
\begin{equation*}
\xi'(1)=\left\langle (-\Delta)^s(u_\epsilon+(c\phi_{1,s}-u_\epsilon)^+), (c\phi_{1,s}-u_\epsilon)^+)\right\rangle-\int_{\Omega}\frac{K_\epsilon}{(c\phi_{1,s}+\epsilon)^\delta}(c\phi_{1,s}-u_\epsilon)^+\,\mathrm{d}x.
\end{equation*}
Using the following well-known convexity inequality
\begin{equation*}
(-\Delta)^s(\phi(u))\leq \phi'(u)(-\Delta)^s u\mbox{ for any convex Lipschitz  function }\phi,
\end{equation*}
we have
\begin{eqnarray*}
\xi'(1)&\leq \left\langle (-\Delta)^s(u_\epsilon+(c\phi_{1,s}-u_\epsilon)), (c\phi_{1,s}-u_\epsilon)^+\right\rangle-\int_{\Omega}\frac{K_\epsilon}{(c\phi_{1,s}+\epsilon)^\delta}(c\phi_{1,s}-u_\epsilon)^+\,\mathrm{d}x\\
&=\left\langle (-\Delta)^s(c\phi_{1,s})-\frac{K_\epsilon}{(c\phi_{1,s}+\epsilon)^\delta}, (c\phi_{1,s}-u_\epsilon)^+)\right\rangle<0
\end{eqnarray*}
if the support of $(c\phi_{1,s}-u_\epsilon)^+$ has non $0$-measure. So we get a contradiction in this case and $c\phi_{1,s}\leq u_\epsilon$. Thus, $E_\epsilon$ is G\^ateaux differentiable at $u_\epsilon$ and $u_\epsilon$ satisfies in the sense of distributions:
\begin{eqnarray*}
\left\{\begin{array}{lll}
&(-\Delta)^s u_\epsilon=\frac{K_\epsilon(x)}{(u_\epsilon+\epsilon)^\delta}\mbox{ in }\Omega,\\
& u_\epsilon>0\mbox{ in }\Omega,\\
& u_\epsilon=0\mbox{ in }\R^N\backslash\Omega.
\end{array}\right.
\end{eqnarray*}
From \cite[Proposition1.1, p.~277]{Ros-oton-serra-JMPA},  we deduce that $u_\epsilon\in C^s(\R^N)$. Now we prove that $u_\epsilon$ is monotone increasing as $\epsilon\downarrow 0^+$ by a comparison argument (that we will use throughout the paper, refered as the comparison principle): Let $0<\epsilon'<\epsilon$. Then, 
\begin{equation*}
(-\Delta)^s(u_{\epsilon'}-u_\epsilon)-\frac{K_{\epsilon'}}{(u_{\epsilon'}+\epsilon')^\delta}+\frac{K_\epsilon}{(u_\epsilon+\epsilon)^\delta}=0.
\end{equation*}
Let $x_0=\displaystyle\arg\min_{\overline{\Omega}}(u_{\epsilon'}-u_\epsilon)$ and assume that $x_0\in \Omega$. Hence, $(u_{\epsilon'}-u_\epsilon)(x_0)\leq 0$. Then,

\begin{eqnarray*}
&(-\Delta)^s(u_{\epsilon'}-u_\epsilon)(x_0)-\frac{K_{\epsilon'}(x_0)}{(u_{\epsilon'}(x_0)+\epsilon')^\delta}+\frac{K_\epsilon(x_0)}{(u_\epsilon(x_0)+\epsilon)^\delta}\leq\\
&C(N,s)\int_{\R^N}\frac{(u_{\epsilon'}-u_\epsilon)(x_0)-(u_{\epsilon'}-u_\epsilon)(y)}{\vert x_0-y\vert^{N+2s}}\,\mathrm{d}y-\frac{K_\epsilon(x_0)}{(u_{\epsilon'}(x_0)+\epsilon')^\delta}+\frac{K_\epsilon(x_0)}{(u_\epsilon(x_0)+\epsilon)^\delta}< 0
\end{eqnarray*}
from which we get a contradiction. Therefore, we get a contradiction and $u_{\epsilon'}> u_\epsilon$ in $\Omega$. Thus, we infer that $u=\displaystyle\lim_{\epsilon\downarrow 0^+}u_\epsilon\geq c\phi_{1,s}$ and satisfies in the sense of distributions
\begin{eqnarray}\label{eq1.2}
\left\{\begin{array}{lll}
&(-\Delta)^s u=\frac{K(x)}{u^\delta}\mbox{ in }\Omega,\\
& u>0\mbox{ in }\Omega,\\
& u=0\mbox{ in }\R^N\backslash\Omega.
\end{array}\right.
\end{eqnarray}
Furthermore, $u\in \tilde{H}^s(\Omega)$. Indeed,
\begin{eqnarray*}
\int_{\R^N}((-\Delta)^{s/2}u_\epsilon)^2\mathrm{d}x&=\int_{\Omega}K_\epsilon\frac{u_\epsilon}{(u_\epsilon+\epsilon)^\delta}\mathrm{d}x\leq\int_{\Omega}K_\epsilon(x)u_\epsilon(x)^{1-\delta}\,\mathrm{d}x\\
&\leq C\left(\int_{\Omega}d(x)^{2(s-\beta-s\delta)}\right)^{1/2}\left(\int_\Omega\left(\frac{u_\epsilon}{d(x)^s}\right)^2\mathrm{d}x\right)^{1/2}
\end{eqnarray*}
from which together with the Hardy inequality it follows
\begin{equation*}
\displaystyle\sup_{\epsilon>0}\int_{\R^N}((-\Delta)^{s/2}u_\epsilon)^2\mathrm{d}x<\infty.
\end{equation*}
Therefore, $u$ being the minimal solution to $(P_s)$, any weak solution to  $(P_s)$ is in $\tilde{H}^s(\Omega)$. Then, using again the comparison principle above, it is easy to prove that $u$ is the unique weak solution to $(P_s)$.
We now prove the upper estimate. Let $G_s(x,y)$ the Green function associated to $(-\Delta)^s$ with homogeneous Dirichlet conditions in $\Omega$. Then we have:
\begin{equation}\label{green-eps}
u_\epsilon=\int_{\Omega}\frac{G_s(x,y)K_\epsilon(y)}{(u_\epsilon(y)+\epsilon)^\delta}\,\mathrm{d}y.
\end{equation}
From \cite[Theorem 1.1, p.~467]{Chen-Song}, we have the following estimates on the Kernel $G_s$:
\begin{equation*}
0\leq G_s(x,y)\leq C\frac{\min\{d^s(x)d^s(y),\vert x-y\vert^sd^s(x)\}}{\vert x-y\vert^N}.
\end{equation*}
Plugging the above inequality into the integral representation \eqref{green-eps}, we obtain for some positive constant $\tilde C$ independent of $\epsilon$:
\begin{equation*}
\frac{u_\epsilon(x)}{d(x)^s}\leq \tilde C\int_\Omega\frac{\min\{d^s(y),\vert x-y\vert^s\}}{\vert x-y\vert^N}\frac{d(y)^{-\beta}}{(u_\epsilon(y)+\epsilon)^\delta}\,\mathrm{d}y.
\end{equation*}
Distinguishing the cases $|x-y|\leq d(y)$ and $d(y)\leq |x-y|$ and noticing that $u_\epsilon\geq c\phi_{1,s}$, we then infer that
\begin{equation*}
\frac{u_\epsilon(x)}{d^s(x)}\leq C\int_\Omega\frac{\min\{d^s(y),\vert x-y\vert^s\}}{\vert x-y\vert^N}d(y)^{-\beta-s\delta}\,\mathrm{d}y\leq C\int_\Omega\vert x-y\vert^{s(1-\delta)-\beta-N}\,\mathrm{d}y<\infty.
\end{equation*}
Thus passing to the limit as $\epsilon\to 0^+$, we infer that 
\begin{equation*}
c\phi_{1,s}\leq u\leq C\phi_{1,s}
\end{equation*}
for some positive constants $c,C$. This completes the proof of the existence and the asymptotic behaviour of the solution to $(P_s)$ in the case $\frac{\beta}{s}+\delta<1$.

Let us consider the case $\frac{\beta}{s}+\delta=1$. Suppose first that $\beta>0$. Following closely the proof of \cite[Proposition 1.2.9, p.~26-30]{Abatangelo}, we can prove the following extension result:

Let $w\in \tilde{H}^s(\Omega)$ be the solution to
\begin{eqnarray}\label{abatangelo-log}
\left\{\begin{array}{ll}
&(-\Delta)^s w=\frac{1}{d(x)^s}\ln^{-\alpha}\left(\frac{A}{d(x)}\right)\mbox{ in }\Omega\\
&w=0\mbox{ in }\R^N\backslash\Omega
\end{array}\right.
\end{eqnarray}
where $0<\alpha<1$ and $A\geq {\rm diam}(\Omega)$. Then, there exists $c_1, c_2>0$ such that
\begin{equation}\label{eq1.3}
c_1d(x)^s\ln^{1-\alpha}\left(\frac{A}{d(x)}\right)\leq w(x)\leq c_2d(x)^s\ln^{1-\alpha}\left(\frac{A}{d(x)}\right),\quad\forall x\in\Omega.
\end{equation}
Furthermore, for any $\alpha_0<1$, $c_1$ and $c_2$ are uniform for $0\leq \alpha\leq \alpha_0$.

Based on this extension result, we prove the existence and behaviour of a unique solution to $(P_s)$. First we prove the existence and uniqueness of the approximated solution $u_\epsilon$. For that we argue as in \cite{peral-al}. Let $C_0(\overline{\Omega})^+$ be the positive cone of $C_0(\overline{\Omega})$ and set the map
\begin{equation*}
T_\epsilon\,:\, C_0(\overline{\Omega})^+\mapsto C_0(\overline{\Omega})^+
\end{equation*}
defined by
\begin{equation*}
T_\epsilon(v)=(-\Delta)^{-s}\left(\frac{K_\epsilon}{(v+\epsilon)^\delta}\right),\quad\mbox{for any }v\in C_0(\overline{\Omega})^+.
\end{equation*}
From \cite[Proposition1.1, p.~277]{Ros-oton-serra-JMPA}, we get for some postive constants $C, C_\epsilon$:
\begin{equation*}
\Vert T_\epsilon v\Vert_{C_0(\overline{\Omega})}\leq \Vert T_\epsilon v\Vert_{C^s(\overline{\Omega})}\leq C\left\Vert \frac{K_\epsilon}{\epsilon^\delta}\right\Vert_{L^\infty(\Omega)}\leq C_\epsilon.
\end{equation*}
Therefore, $T_\epsilon$ is a continuous and compact operator from $M_\epsilon=\left\{\phi\in C_0(\overline{\Omega})^+\,|\, \Vert\phi\Vert_{C_0(\overline{\Omega})}\leq C_\epsilon\right\}$ onto $M_\epsilon$. Then, applying the Schauder fix point Theorem, we infer the existence of $u_\epsilon$, solution to $(P_\epsilon)$. The uniqueness of $u_\epsilon$ follows from the same arguments used in the case $\frac{\beta}{s}+\delta<1$. As in the case $\frac{\beta}{s}+\delta<1$, we have also for some positive constant $c$ independent of $\epsilon$:
\begin{equation*}
c\phi_{1,s}\leq u_\epsilon\quad\mbox{and}\quad c\phi_{1,s}\leq u=\displaystyle\lim_{\epsilon\downarrow 0^+}u_\epsilon
\end{equation*}
and $u$ satisfies  \eqref{eq1.2} in the sense of distributions. Next we establish the asymptotic behaviour of $u$ near the boundary. Precisely, we aim to show that for some constant $D>0$ large enough,
\begin{equation}\label{eq1.4}
\frac{1}{D}\phi_{1,s}\ln^{\frac{1}{\delta+1}}\left(\frac{2}{\phi_{1,s}}\right)\leq u\leq D\phi_{1,s}\ln^{\frac{1}{\delta+1}}\left(\frac{2}{\phi_{1,s}}\right).
\end{equation}
We can assume that $\delta>0$ (if $\delta=0$, \eqref{eq1.4} follows from \eqref{eq1.3} with $\alpha=0$). For that we iterate some estimates from \eqref{eq1.3} for suitable values of $\alpha\in (0,1)$ and $\alpha_0$. Precisely, from \eqref{eq1.3} and \eqref{eq1.1}, the following estimates hold for some positive constants $M, C_0$ large enough and $0\leq \alpha\leq \frac{1+\delta^3}{1+\delta}<1$:
\begin{equation}\label{eq1.5}
\frac{1}{M}\phi_{1,s}\ln^{1-\alpha}\left(\frac{2}{\phi_{1,s}}\right)\leq (-\Delta)^{-s}\left(\frac{1}{\phi_{1,s}}\ln^{-\alpha}\left(\frac{2}{\phi_{1,s}}\right)\right)\leq M\phi_{1,s}\ln^{1-\alpha}\left(\frac{2}{\phi_{1,s}}\right),
\end{equation}
\begin{equation}\label{eq1.6}
c\phi_{1,s}\leq u,\quad \frac{1}{C_0}\phi_{1,s}^{\frac{-\beta}{s}}\leq K(x)\leq C_0\phi_{1,s}^{\frac{-\beta}{s}}.
\end{equation}
From \eqref{eq1.6}, we get that for any $\epsilon>0$
\begin{equation*}
(-\Delta)^s u_\epsilon\leq \frac{C_0\phi_{1,s}^{\frac{-\beta}{s}}}{(c\phi_{1,s})^\delta}= C_0c^{-\delta}\phi_{1,s}^{-1}.
\end{equation*}
Using \eqref{eq1.3} and the comparison principle, we infer that 
\begin{equation}\label{eq1.7}
u_\epsilon\leq C_0c^{-\delta}M\phi_{1,s}\ln\left(\frac{2}{\phi_{1,s}}\right)\mbox{and by taking }\epsilon\to0^+\;u\leq C_0c^{-\delta}M\phi_{1,s}\ln\left(\frac{2}{\phi_{1,s}}\right)
\end{equation}
from which we obtain $u\in C_0(\overline{\Omega})$. Again using the equation \eqref{eq1.2} satisfied by $u$ and plugging the estimate \eqref{eq1.7}, we obtain:
\begin{equation*}
(-\Delta)^s u\geq \frac{\frac{1}{C_0}\phi_{1,s}^{-\frac{\beta}{s}}}{\left(C_0c^{-\delta}M\phi_{1,s}\ln\left(\frac{2}{\phi_{1,s}}\right)\right)^\delta}=C_0^{-1-\delta}c^{\delta^2}M^{-\delta}\phi_{1,s}^{-1}\ln^{-\delta}\left(\frac{2}{\phi_{1,s}}\right)
\end{equation*}
from which together with the comparison principle it follows that
\begin{equation*}
u\geq \frac{c^{\delta^2}}{(MC_0)^{\delta +1}}\phi_{1,s}\ln^{1-\delta}\left(\frac{2}{\phi_{1,s}}\right).
\end{equation*}
Iterating these estimates, we get for any $p\in\N$:
\begin{equation}\label{esti-iteration}
\frac{c^{\delta^{2p+2}}}{(MC_0)^{1+\delta+\cdot\cdot+\delta^{2p+1}}}\phi_{1,s}\ln^{1-\delta+\delta^2+\cdot\cdot-\delta^{2p+1}}\left(\frac{2}{\phi_{1,s}}\right)\leq u\leq 
\frac{(MC_0)^{1+\delta+\cdot\cdot+\delta^{2p+2}}}{c^{\delta^{2p+3}}}\phi_{1,s}\ln^{1-\delta+\delta^2+\cdot\cdot+\delta^{2p+2}}\left(\frac{2}{\phi_{1,s}}\right).
\end{equation}
Passing to the limit as $p\to\infty$, we obtain that
\begin{equation*}
\left(\frac{1}{MC_0}\right)^{\frac{1}{1-\delta}}\phi_{1,s}\ln^{\frac{1}{\delta+1}}\left(\frac{2}{\phi_{1,s}}\right)\leq u\leq (MC_0)^{\frac{1}{1-\delta}}\phi_{1,s}\ln^{\frac{1}{\delta+1}}\left(\frac{2}{\phi_{1,s}}\right).
\end{equation*}
Let us finally consider the case $\beta=0$ and $\delta=1$. In this case $u$ satisfies 
\begin{eqnarray*}
\left\{\begin{array}{lll}
&(-\Delta)^s u=\frac{1}{u}\mbox{ in }\Omega,\\
& u>0\mbox{ in }\Omega,\\
& u=0\mbox{ in }\R^N\backslash\Omega.
\end{array}\right.
\end{eqnarray*}
From \eqref{eq1.3}, we have for $M_0>0$ large enough,
\begin{equation*}
\frac{1}{M_0}\phi_{1,s}\ln^{\frac{1}{2}}\left(\frac{2}{\phi_{1,s}}\right)\leq w_0\eqdef(-\Delta)^{-s}\left(\frac{1}{\phi_{1,s}\ln^{\frac{1}{2}}\left(\frac{2}{\phi_{1,s}}\right)}\right)\leq M_0\phi_{1,s}\ln^{\frac{1}{2}}\left(\frac{2}{\phi_{1,s}}\right).
\end{equation*}
Then for a positive constant $C$ large enough, we have
\begin{equation*}
(-\Delta)^s(Cw_0)\geq \frac{1}{Cw_0}\;\mbox{ and }\; (-\Delta)^s\left(\frac{1}{C}w_0\right)\leq \frac{C}{w_0}\;\mbox{ in }\Omega.
\end{equation*}
Therefore again by the comparison principle, $\frac{w_0}{C}\leq u\leq Cw_0$. Thus for a constant $C_1$ large enough,
\begin{equation*}
\frac{1}{C_1}\phi_{1,s}\ln^{\frac{1}{2}}\left(\frac{2}{\phi_{1,s}}\right)\leq u\leq C_1\phi_{1,s}\ln^{\frac{1}{2}}\left(\frac{2}{\phi_{1,s}}\right).
\end{equation*}
Next, we consider the case $\frac{\beta}{s}+\delta>1$ which is equivalent to $\beta>(1-\delta)s$. Let $w\in L^1(\Omega)\cap C^s_{{\rm loc}}(\Omega)$ be the function satisfying:
\begin{equation*}
 (-\Delta)^s w=\frac{1}{d(x)^{\alpha\delta+\beta }}=\frac{1}{d(x)^{2s-\alpha }}\;\mbox{ with }\alpha=\frac{2s-\beta}{\delta +1}.
\end{equation*}
Then, from\cite[Proposition 1.2.9]{Abatangelo}, there exists $M_1>0$ such that 
\begin{equation*}
\frac{1}{M_1}\phi_{1,s}^{\frac{\alpha}{s}}\leq w\leq M_1\phi_{1,s}^{\frac{\alpha}{s}}
\end{equation*}
and for a constant $C>0$ large enough
\begin{equation*}
(-\Delta)^s(Cw)\geq \frac{K(x)}{(Cw)^\delta}\;\mbox{ and }\; (-\Delta)^s\left(\frac{w}{C}\right)\leq \frac{C^\delta K(x)}{w^\delta}.
\end{equation*}
Then, we conclude as above by the comparison principle that
\begin{equation*}
\frac{w}{C}\leq u\leq C w\;\mbox{ thus }\; \frac{1}{CM_1}\phi_{1,s}^{\frac{2s-\beta}{s(\delta +1)}}\leq u\leq CM_1 \phi_{1,s}^{\frac{2s-\beta}{s(\delta +1)}}.
\end{equation*}
Finally, observe that $u\in \tilde{H}^s(\Omega)$ if and only if $\int_\Omega K(x) u^{1-\delta}\,\mathrm{d}x$. Due to the behaviour of $u$ with respect to $d$, it reduces to the necessary and sufficient condition $2\beta+\delta(2s-1)<1+2s$. This completes the proof of assertions (i)-(ii). Assertion (iii) follows from the last assertion in (16) of  \cite[Proposition 1.2.9, p.~7]{Abatangelo}. This completes the proof of Theorem~\ref{sing-prob}.
\qed

We now prove Theorem~\ref{sing-prob2}. To this aim, we recall the following regularity results from \cite[Theorem 2.9, p.~79]{Silvestre-CPAM} (see also \cite[Proposition 2.3, p. 280]{Ros-oton-serra-JMPA}):
\begin{prop}\label{regul-holder}
Let $w\in C^\infty(\R^N)$, $h\in L^\infty(\Omega)$ such that
\begin{equation*}
(-\Delta)^sw=h\mbox{ in }B_1.
\end{equation*}
Then, $\forall\tilde{\beta}\in (0,2s)$
\begin{equation*}
\Vert w\Vert_{C^{\tilde{\beta}}(\overline{B}_{1/2})}\leq C\left(\Vert w\Vert_{L^\infty(\R^N)}+\Vert h\Vert_{L^\infty(B_1)}\right).
\end{equation*}
\end{prop}
As in \cite[Corollary 2.5, p.~280]{Ros-oton-serra-JMPA}, we have also
\begin{coro}\label{coro-holder}
Assume that $w\in C^\infty(\R^N)$ is a solution to
\begin{equation*}
(-\Delta)^s w=h\mbox{ in } B_1.
\end{equation*}
Then, for every $\tilde{\beta}\in (0,2s)$,
\begin{equation*}
\Vert w\Vert_{C^{\tilde{\beta}}(\overline{B}_{1/2})}\leq C\left(\Vert(1+\vert x\vert)^{-N-2s}w\Vert_{L^1(\R^N)}+\Vert w\Vert_{C^{\tilde{\beta}}(\overline{B}_2)}+\Vert h\Vert_{C^{\tilde{\beta}}(\overline{B}_2)}\right).
\end{equation*}
where the constant $C$ depends only on $n$, $s$ and $\tilde{\beta}$.
\end{coro}
Using the regularity estimates above as in \cite[Lemma 2.9, p.~281]{Ros-oton-serra-JMPA}, we establish the following result:
\begin{lemm}\label{regu-esti-local}
Let $x_0\in \Omega$ and $R=\frac{d(x_0)}{2}$. Then, for any $\tilde{\beta}\in (0,2s)$, there exists $C=C(\tilde{\beta},s,\Omega)>0$ and $C_\epsilon=C_\epsilon(\tilde{\beta},s,\Omega,\epsilon)$ such that:
\begin{itemize}
\item[i)] If $\frac{\beta}{s}+\delta<1$, $
\Vert u\Vert_{C^{\tilde{\beta}}(\overline{B_R(x_0)})}\leq CR^{s-\tilde{\beta}}$;
\item[ii)] if $\frac{\beta}{s}+\delta=1$, $\Vert u\Vert_{C^{\tilde{\beta}}(\overline{B_R(x_0)})}\leq C_\epsilon R^{s-\tilde{\beta}-\epsilon}, \;\forall\epsilon>0\; \mbox{ small enough}$;
\item[iii)] If $\frac{\beta}{s}+\delta>1$, $
\Vert u\Vert_{C^{\tilde{\beta}}(\overline{B_R(x_0)})}\leq CR^{\frac{2s-\beta}{\delta+1}-\tilde{\beta}}$.
\end{itemize}
\end{lemm}
{\bf Proof of Lemma~\ref{regu-esti-local}.}\\
Without loss of generality, we may assume $u\in C^\infty(\R^N)$. Indeed, we can regularize by the usual mollification technique: $u*\eta_\epsilon$ where $\eta_\epsilon$ is the standard mollifier.

Note that $B_R(x_0)\subset B_{2R}(x_0)\subset \Omega$. As in the proof of Lemma 2.9 p.~282 of \cite{Ros-oton-serra-JMPA},  we define the rescaled function $\tilde{u}$:
\begin{equation*}
\tilde{u}(y)=u(x_0+Ry)
\end{equation*}
and let $g(y)=\frac{K(y)}{u^\delta(y)}$. Then, we have 
\begin{equation*}
(-\Delta)^s\tilde{u}(y)=R^{2s}g(x_0+Ry)\;\mbox{ in }\; B_1.
\end{equation*}
Thus, using Corollary~\ref{coro-holder}, we obtain:
\begin{equation*}
\Vert \tilde{u}\Vert_{C^{\tilde{\beta}}(\overline{B_{1/4}})}\leq C\left(\Vert R^{2s}g(x_0+Ry)\Vert_{L^\infty(B_1)}+\Vert\tilde{u}\Vert_{L^\infty(B_1)}+\Vert(1+\vert x\vert)^{-N-2s}\tilde{u}\Vert_{L^1(\R^N)}\right).
\end{equation*}
We distinguish the following cases:\\
{\it Case 1:} $\frac{\beta}{s}+\delta<1$. From Theorem~\ref{sing-prob}, we have $\tilde{u}(x)\leq CR^s\mbox{ for any }x\in B_1$. Then,
\begin{equation*}
\Vert \tilde{u}\Vert_{L^\infty(B_1)}\leq CR^s\;\mbox{ and }\;\Vert R^{2s}g(x_0+Ry)\Vert_{L^\infty(B_1)}\leq CR^{2s-\beta-\delta s}.
\end{equation*}
Furthermore observing that
\begin{equation*}
\tilde{u}(y)\leq CR^s(1+\vert y\vert^s)\quad\forall y\in \R^N,
\end{equation*}
we infer that
\begin{equation*}
\Vert(1+\vert y\vert)^{-N-2s}\tilde{u}\Vert_{L^1(\R^N)}\leq CR^s.
\end{equation*}
Therefore, gathering the above estimates and since $\beta+\delta s<s$, we obtain that $\Vert\tilde{u}\Vert_{C^{\tilde{\beta}}(\overline{B_{1/4}})}\leq CR^s$ with $C$ independent of $R$. Then, we conclude that
\begin{equation*}
\Vert\tilde{u}\Vert_{C^{\tilde{\beta}}(\overline{B_{R/4}(x_0)})}\leq CR^{s-\tilde{\beta}}.
\end{equation*}
This completes the proof  of assertion (i) of Lemma~\ref{regu-esti-local}.\\
{\it Case 2:} $\frac{\beta}{s}+\delta=1$. In this case, we recall that from Theorem~\ref{sing-prob} $\forall\epsilon>0$ small enough, $u\leq C_\epsilon d^{s-\epsilon}$ and then 
\begin{equation*}
\forall\epsilon>0\,\mbox{(small enough) }\,,\;\Vert\tilde{u}\Vert_{L^\infty(B_1)}\leq C_\epsilon R^{s-\epsilon},
\end{equation*}
\begin{equation*}
\Vert R^{2s} g(x_0+Ry)\Vert_{L^\infty(B_1)}\leq C_\epsilon R^{-\beta-\delta(s-\epsilon)+2s}\mbox{ and}
\end{equation*}
\begin{equation*}
\Vert (1+\vert y\vert)^{-N-2s}\tilde{u}\Vert_{L^1(\R^N)}\leq C_\epsilon R^{s-\epsilon}.
\end{equation*}
Therefore, 
\begin{equation*}
\Vert\tilde{u}\Vert_{C^{\tilde{\beta}}(\overline{B_{1/4}})}\leq C_\epsilon R^{s-\epsilon}\,\mbox{ and then }\,\Vert u\Vert_{C^{\tilde{\beta}}(\overline{B_{R/4}(x_0)})}\leq C_\epsilon R^{s-\epsilon-\tilde{\beta}}.
\end{equation*}
Similarly, we deal with {\it Case 3}: $\frac{\beta}{s}+\delta>1$. In this case, from Theorem~\ref{sing-prob}, we have 
\begin{eqnarray*}
&\Vert \tilde u\Vert_{L^\infty(B_1)}\leq CR^{\frac{2s-\beta}{\delta+1}},\; \Vert R^{2s}g(x_0+Ry)\Vert_{L^\infty(B_1)}\leq CR^{2s-\beta-\frac{\delta(2s-\beta)}{\delta+1}}=CR^{\frac{2s-\beta}{\delta+1}}\;\mbox{and }\\
&\Vert (1+\vert x\vert)^{-N-2s}\tilde{u}}\Vert_{L^1(\R^N)}\leq CR^{\frac{2s-\beta}{\delta+1}.
\end{eqnarray*}
Therefore,
\begin{equation*}
\Vert\tilde{u}\Vert_{C^{\tilde{\beta}}(\overline{B_{1/4}})}\leq C R^{\frac{2s-\beta}{\delta+1}}\,\mbox{ and then }\,\Vert u\Vert_{C^{\tilde{\beta}}(\overline{B_{R/4}(x_0)})}\leq CR^{\frac{2s-\beta}{\delta+1}-\tilde{\beta}}.
\end{equation*}
Finally, by a covering argument, the H\"older estimates in Lemma~\ref{regu-esti-local} hold. This completes the proof of the Lemma.
\qed

\noindent{\bf Proof of Theorem~\ref{sing-prob2}.}\\
We again distinguish cases with respect to $\beta$ and $\delta$ and follow the proof of \cite[Proposition 1.1, p~282]{Ros-oton-serra-JMPA}:\\
{\it Case 1: $\frac{\beta}{s}+\delta<1$}. Doing $\tilde{\beta}=s$, we have for $x_0\in \Omega$, 
\begin{equation}\label{eq1.9}
\Vert u\Vert_{C^s(\overline{B_R(x_0)})}\leq C,\mbox{ with }R=\frac{d(x_0)}{2}.
\end{equation}
To get the estimate on all $\Omega$ (and then on all $\R^N$ since $u\equiv 0$ on $\R^N\backslash\Omega$), it is sufficient from \eqref{eq1.9} and interior regularity that follows from \cite[Proposition 1.1]{Ros-oton-serra-JMPA},  to extend \eqref{eq1.9} on $\displaystyle\cup_{x_0\in \Omega_\eta}\overline{B_{2R}(x_0)}\backslash B_R(x_0)$ where $\eta>0$ small enough and $\Omega_\eta=\{x\in\Omega\,:\,d(x)<\eta\}$.

In this regard, let $x,y\in \Omega_\eta$ with $\vert x-y\vert\geq \max(d(x)/2,d(y)/2)$. Suppose that $\frac{\beta}{s}+\delta<1$. Then, for a constant $C>0$ large enough,
\begin{equation*}
\frac{\vert u(x)-u(y)\vert}{\vert x-y\vert^s}\leq \frac{\vert u(x)\vert}{\vert x-y\vert^s}+\frac{\vert u(y)\vert}{\vert x-y\vert^s}\leq 2^s\left(\frac{u(x)}{d(x)^s}+\frac{u(y)}{d(y)^s}\right)\leq C.
\end{equation*}
In the cases $\frac{\beta}{s}+\delta=1$ and $\frac{\beta}{s}+\delta>1$, we obtain respectively
\begin{equation*}
\frac{\vert u(x)-u(y)\vert}{\vert x-y\vert^{s-\epsilon}}\leq C_\epsilon\mbox{ and }\frac{\vert u(x)-u(y)\vert}{\vert x-y\vert^{\frac{2s-\beta}{\delta+1}}}\leq C.
\end{equation*}
This completes the proof of Theorem~\ref{sing-prob2}.
\qed
\section{Global bifurcation results}\label{section2}
We now prove Theorems~\ref{main-bifurcation} and \ref{bifurcation-example}. We start by the following result which states the existence and the regularity of the branch of minimal (classical) solutions to $(P_\lambda)$ for $\lambda\in (0,\Lambda)$ with $\Lambda>0$.
\begin{prop}\label{minimal_branch}
Assume {\bf (f1)}-{\bf (f3)}. For any $\delta\geq 0$ and $0\leq \beta<2s$, there exists $\Lambda>0$ such that
\begin{itemize}
\item[i)] For any $0<\lambda<\Lambda$, $(P_\lambda)$ admits a minimal solution $u_\lambda\in C_{\phi_{\delta,\beta}}^+(\Omega)$. For $\lambda>\Lambda$, here is no weak bounded solution. 
\item[ii)] There exists $C>0$ such that for $\lambda>0$ small enough, $u_\lambda$ is the unique solution in ${\mathcal S}\cap\R^+\times\{u\in L^\infty(\Omega)\,:\, \Vert u\Vert_{L^\infty(\Omega)}\leq C\}$.
\item[iii)]  $u_\lambda\in \tilde{H}^s(\Omega)$ if and only if $2\beta+\delta(2s-1)<1+2s$. Assuming $2\beta+\delta(2s-1)<1+2s$, then $(0,\Lambda)\ni \lambda\mapsto u_\lambda\in C_0(\overline{\Omega})$ is of class $C^2$.
Furthermore, $u_\lambda\to u_\Lambda$ in $\tilde{H}^s(\Omega)$ as $\lambda\to\Lambda^-$, where $u_\Lambda$ is a weak solution to $(P_\lambda)$ for $\lambda=\Lambda$. 
\end{itemize}
\end{prop}
{\bf Proof of Proposition~\ref{minimal_branch}.}\\
We start by showing the existence of $u_\lambda$ for $\lambda>0$ small. Let ${\underline{u}}_\lambda$ be the solution to
\begin{eqnarray}\label{strictsub}
\displaystyle\left\{\begin{array}{lll}
& (-\Delta)^su=\frac{\lambda K(x)}{u^\delta}\mbox{ in }\Omega,\\
& u>0\mbox{ in }\Omega,\\
& u=0\mbox{ in }\R^N\backslash\Omega
\end{array}\right.
\end{eqnarray}
${\underline{u}}_\lambda=\lambda^{\frac{1}{\delta+1}} {\underline{u}}_1$ where ${\underline{u}}_1$ is given by Theorem~\ref{sing-prob}. Clearly, ${\underline{u}}_\lambda$ is a strict subsolution to $(P_\lambda)$. Let $U\in C^s(\R^N)\cap \tilde{H}^s(\Omega)$ be the unique solution to 
\begin{eqnarray}\label{supersol}
\displaystyle\left\{\begin{array}{lll}
& (-\Delta)^su=1\mbox{ in }\Omega,\\
& u>0\mbox{ in }\Omega,\\
& u=0\mbox{ in }\R^N\backslash\Omega
\end{array}\right.
\end{eqnarray}
Then, setting $\bar{u}_\lambda\eqdef {\underline{u}}_\lambda+MU$ for some $M>1$ and  letting $\lambda_0>0$, we have
\begin{equation*}
(-\Delta)^s(\bar{u}_\lambda)=\frac{\lambda K(x)}{{\underline{u}}_\lambda^\delta}+M\geq \frac{\lambda K(x)}{({\underline{u}}_\lambda+MU)^\delta}+\lambda f({\underline{u}}_\lambda+MU)
\end{equation*}
if 
\begin{equation}\label{eqM}
\frac{M}{\displaystyle\max_{\lambda\leq\lambda_0}\Vert f({\underline{u}}_\lambda+MU)\Vert_{L^\infty(\Omega)}}\geq \lambda_0\mbox{ and } \lambda\leq\lambda_0.
\end{equation}
Taking $\lambda_0>0$ small enough such that \eqref{eqM} is verified, $\bar{u}_\lambda$ is then a supersolution to $(P_\lambda)$ for $\lambda\leq \lambda_0$. Next, we define the following iterative scheme ($n\geq 1$):
\begin{eqnarray*}
\displaystyle\left\{\begin{array}{lll}
& (-\Delta)^su_n+\lambda Cu_n-\frac{\lambda K(x)}{u_n^\delta}=\lambda Cu_{n-1}+\lambda f(u_{n-1})\mbox{ in }\Omega,\\
& u>0\mbox{ in }\Omega,\\
& u=0\mbox{ in }\R^N\backslash\Omega
\end{array}\right.
\end{eqnarray*}
with $u_0={\underline{u}_\lambda}$ and $C=C(\lambda_0)>0$ large enough such that $t\to Ct+ f(t)$ is increasing on $[0,\Vert \bar{u}_{\lambda_0}\Vert_{L^\infty(\Omega)}]$.

Using the comparison principle as above (for that we remark that the operator $(-\Delta)^su+\lambda Cu-\frac{\lambda K(x)}{u^\delta}$ is monotone) and adapting slightly the proof of Theorem~\ref{sing-prob} to get the existence of $u_n$, it is easy to prove that $(u_n)_{n\geq 0}\subset C^s(\R^N)\cap C_{\phi_{\delta,\beta}}^+(\Omega)$ and is increasing. Furthermore, for any $0<\lambda\leq \lambda_0$, ${\underline{u}_\lambda}\leq u_n\leq {\underline{u}_\lambda}+MU$. From Theorems~\ref{sing-prob} and \ref{sing-prob2}, $\displaystyle\sup_{n\in\N}\Vert u_n\Vert_{C^{\gamma}(\R^N)}\leq C_0$ for some $\gamma=\gamma(\delta, \beta,\lambda_0)$ and $C_0=C_0(\delta, \beta,\lambda_0)$ large enough. Therefore,
\begin{equation*}
u_n\to u\mbox{ in }C(\R^N)\mbox{ as }n\to\infty\mbox{ and $u$ satisfies}
\end{equation*}
\begin{equation*}
(-\Delta)^s u=\lambda\left(\frac{K(x)}{u^\delta}+f(u)\right)
\end{equation*}
in the sense of distributions. Now, we set 
\begin{equation*}
\Lambda=\displaystyle\sup\left\{\lambda>0\,:\,\exists\mbox{ a weak bounded solution to }(P_\lambda)\right\}.
\end{equation*}
Obviously, from above we have $\Lambda>0$. Furthermore, from assumption {\bf (f3)}, $\Lambda<\infty$. In addition, for any $0<\lambda<\Lambda$, there exists $u_\lambda\in C_{\phi_{\delta,\beta}}^+(\Omega)$ a minimal solution to $(P_\lambda)$ (for any $\lambda\in (0,\Lambda)$, take ${\underline{u}}_\lambda$ as a subsolution and $v_{\lambda'}$, solution to $(P_{\lambda'})$, with suitable $\lambda<\lambda'<\Lambda$ as a supersolution). This completes the proof of assertion i). Let us prove assertion ii). Note that from the comparison principle, $(0,\Lambda)\ni\lambda\mapsto u_\lambda\in C_0(\overline{\Omega})$ is increasing and $u_\lambda\to 0$ in $C_0(\overline{\Omega})$ as $\lambda\to 0^+$. Since that for $\lambda_0>0$ small enough and uniformly with respect to $x\in\overline{\Omega}$,
\begin{equation}\label{uniqueness-small}
t\to  K(x) t^{-\delta} + f(t)\mbox{ is decreasing for } t\leq C_0\eqdef\Vert u_{\lambda_0}\Vert_{L^\infty(\Omega)}
\end{equation}
and from comparison principle, we get that $u_\lambda$ is the unique solution to $(P_\lambda)$ in $\left\{u\in C_0(\overline{\Omega})\,:\,\Vert u\Vert_{L^\infty(\Omega)}\leq C_0\right\}$. Indeed, consider $v_\lambda$ solution to $(P_\lambda)$ satisfying $\Vert v_\lambda\Vert_{L\infty(\Omega)}\leq C_0$. Let $x_0\in \Omega$ satisfying $u_\lambda(x_0)-v_\lambda(x_0)=\displaystyle\min_{\Omega} u_\lambda-v_\lambda$. Suppose that $u_\lambda(x_0)-v_\lambda(x_0)<0$. Then, since $u_\lambda$ is the minimal solution,
\begin{equation*}
(-\Delta)^s(u_\lambda-v_\lambda)(x_0)=\int_{\Omega}\frac{(u_\lambda-v_\lambda)(x_0)-(u_\lambda-v_\lambda)(y)}{\vert x_0-y\vert^{N+2s}}\,\mathrm{d}y<0.
\end{equation*}
On the other hand, from \eqref{uniqueness-small}
\begin{equation*}
\frac{\lambda K(x_0)}{u_\lambda^\delta(x_0)}+\lambda f(u_\lambda(x_0))-\frac{\lambda K(x_0)}{v_\lambda^\delta(x_0)}-\lambda f(v_\lambda(x_0))>0
\end{equation*}
from which we get a contradiction. This yields $v_\lambda\equiv u_\lambda$.

Finally let us prove iii). First, we observe that 
\begin{equation*}
2\beta+\delta(2s-1)<1+2s\Leftrightarrow u_\lambda\in \tilde{H}^s(\Omega),\,\forall\lambda\in (0,\Lambda)\Leftrightarrow{\underline{u}_\lambda}\in \tilde{H}^s(\Omega)\;\forall\lambda>0.
\end{equation*}
Indeed, 
\begin{equation*}
\int_{\Omega}K(x)u_\lambda^{1-\delta}\,\mathrm{d}x<\infty\Leftrightarrow\int_{\Omega}K(x){\underline{u}}_\lambda^{1-\delta}\,\mathrm{d}x<\infty\Leftrightarrow 2\beta+\delta(2s-1)<1+2s.
\end{equation*}
Assuming $2\beta+\delta(2s-1)<1+2s$, we consider the eigenvalue problem:
\begin{equation*}
\Lambda_1(\lambda)=\displaystyle\inf_{\phi\in \tilde{H}^s(\Omega),\,\int_{\Omega}\phi^2\,\mathrm{d}x=1}\left\{\int_{\mathbb{R}^N}\vert(-\Delta)^{\frac{s}{2}}\phi\vert^2\,\mathrm{d}x+\lambda\delta\int_{\Omega}\frac{K(x)\phi^2}{u_\lambda^{\delta+1}}\,\mathrm{d}x-\lambda\int_{\Omega}f'(u_\lambda)\phi^2\,\mathrm{d}x\right\}.
\end{equation*}
From the Hardy inequality which implies that any $u\in \tilde{H}^s(\Omega)$ satisfies $\int_{\Omega}\frac{u^2}{d(x)^{2s}}<\infty$), and the compact embedding of $\tilde{H}^s(\Omega)$ in $L^2(\Omega)$, $\Lambda_1(\lambda)$ is achieved on some $\phi_\lambda\in \tilde{H}^s(\Omega)$ with $\int_{\Omega}\phi_\lambda^2\,\mathrm{d}x=1$. Furthermore, $\phi_\lambda$ satisfies
\begin{equation*}
(-\Delta)^s\phi_\lambda+\frac{\lambda\delta K(x)}{u_\lambda^{\delta+1}}\phi_\lambda=\lambda f'(u_\lambda)\phi_\lambda+\Lambda_1(\lambda)\phi_\lambda
\end{equation*}
that is $\phi_\lambda$ is an eigenfunction associated to the first eigenvalue $\Lambda_1(\lambda)$ of the operator $(-\Delta)^s+\frac{\lambda\delta K(x)}{u_\lambda^{\delta+1}}-\lambda f'(u_\lambda)$. 

We now show that $\Lambda_1(\lambda)$ is a principal eigenvalue, i.e. $\phi_\lambda$ (and any eigenfunction associated to $\Lambda_1(\lambda)$) does not change sign. Assume by contradiction that $\phi_\lambda^+\not\equiv 0$ and $\phi_\lambda^-\not\equiv 0$. Then, 
\begin{eqnarray*}
\begin{array}{llll}
&\left\langle(-\Delta)^s\phi_\lambda,\phi_\lambda^+\right\rangle=C(N,s)\int_{\R^N}\int_{\R^N}\frac{(\phi_\lambda(x)-\phi_\lambda(y))(\phi_\lambda^+(x)-\phi_\lambda^+(y))}{\vert x-y\vert^{N+2s}}\,\mathrm{d}x\mathrm{d}y=\\
&C(N,s)\int_{\R^N}\int_{\R^N}\frac{(\phi_\lambda^+(x)-\phi_\lambda^+(y))(\phi_\lambda^+(x)-\phi_\lambda^+(y))}{\vert x-y\vert^{N+2s}}\,\mathrm{d}x\mathrm{d}y+C(N,s)\int_{\R^N}\int_{\R^N}\frac{\phi_\lambda^-(x)\phi_\lambda^+(y)}{\vert x-y\vert^{N+2s}}\,\mathrm{d}x\mathrm{d}y+\\
&C(N,s)\int_{\R^N}\int_{\R^N}\frac{\phi_\lambda^+(x)\phi_\lambda^-(y)}{\vert x-y\vert^{N+2s}}\,\mathrm{d}x\mathrm{d}y
=-\lambda\delta\int_{\Omega}\frac{K(x)\phi_\lambda^+(x)^2}{u_\lambda(x)^{\delta+1}}\,\mathrm{d}x+\lambda\int_{\Omega}f'(u_\lambda)\phi_\lambda^+(x)^2\,\mathrm{d}x+\Lambda_1(\lambda)\int_\Omega\phi_\lambda^+(x)^2\,\mathrm{d}x.
\end{array}
\end{eqnarray*}
If $\phi_\lambda^-\not\equiv 0$ and $\phi_\lambda^+\not\equiv 0$, it implies that
\begin{equation*}
\int_{\R^N}((-\Delta)^{\frac{s}{2}}\phi_\lambda^+(x))^2\,\mathrm{d}x+\lambda\delta\int_{\Omega}\frac{K(x)\phi_\lambda^+(x)^2}{u_\lambda(x)^{\delta+1}}\,\mathrm{d}x-\lambda\int_{\Omega}f'(u_\lambda)\phi_\lambda^+(x)^2\,\mathrm{d}x<\Lambda_1(\lambda)\int_\Omega\phi_\lambda^+(x)^2\,\mathrm{d}x.
\end{equation*}
This contradicts the definition of $\Lambda_1(\lambda)$. Therefore, $\Lambda_1(\lambda)$ is a principal eigenvalue.  Without loss of generality, we assume that $\phi_\lambda$ is nonnegative. Let us show that $\phi_\lambda$ is positive in $\Omega$. Following the arguments in the proof of \cite[Theorem 3.2, p.~379]{Franzina-Palatucci}, we deduce that $\phi_\lambda\in L^\infty(\Omega)$ and from the local regularity theory that $\phi_\lambda\in C^s_{\rm loc}(\Omega)$. Then, assume by contradiction that $\exists x_0\in \Omega$ such that $\phi_\lambda(x_0)=0$. It follows that
\begin{equation*}
0> 2C(N,s)\int_{\R^N}\frac{\phi_\lambda(x_0)-\phi_\lambda(y)}{\vert x-y\vert^{N+2s}}\,\mathrm{d}y=-\frac{\lambda\delta K(x_0)\phi_\lambda(x_0)}{u_\lambda^{\delta +1}(x_0)}+\lambda f'(u_\lambda(x_0))\phi_\lambda(x_0)+\Lambda_1(\lambda)\phi_\lambda(x_0)=0
\end{equation*}
from which we get a contradiction. Thus $\phi_\lambda>0$ in $\Omega$. Based on this result, we can follow the proof of \cite[Theorem 4.2, p.~382]{Franzina-Palatucci} using the convex inequality stated in \cite[Theorem 4.1, p.~381]{Franzina-Palatucci}, we prove also that $\Lambda_1(\lambda)$ is simple.

From a classical subsolution and supersolution argument, we obtain that $\Lambda_1(\lambda)\geq 0$ for any $\lambda\in (0,\Lambda)$. Indeed, assuming that for some $\lambda\in (0,\Lambda)$, $\Lambda_1(\lambda)<0$. Then, we can show that $u_\lambda-\epsilon\phi_\lambda$ is a strict supersolution to $(P_\lambda)$. Hence we can prove by using the iterative scheme above that there exists a  weak bounded solution to $(P_\lambda)$, $v_\lambda$, such that ${\underline{u}}_\lambda\leq v_\lambda\leq u_\lambda-\epsilon\phi_\lambda$ that contradicts that $u_\lambda$ is a minimal solution.

Next, from the strict convexity of $t\to \lambda(t^{-\delta}+f(t))$, we obtain that $\Lambda_1(\lambda)$ is strictly monotone with respect to $\lambda\in (0,\Lambda)$. Thus, $\Lambda_1(\lambda)>0$ for all $\lambda\in (0,\Lambda)$. Let $g(x,u)\eqdef\left(\frac{K(x)}{u^\delta}+f(u)\right)$. Therefore, from the second statement in {\bf(f3)} and noticing that $\displaystyle\sup_{\lambda\in[\Lambda/2,\Lambda]}\int_{\Omega}K(x)u_\lambda^{1-\delta}<\infty$ when $2\beta+\delta(2s-1)<1+2s$, we have that
\begin{equation*}
\int_{\R^N}((-\Delta)^{\frac{s}{2}}u_\lambda)^2\,\mathrm{d}x=\lambda\int_{\Omega}g(x,u_\lambda)u_\lambda\,\mathrm{d}x\mbox{ and }\int_{\R^N}((-\Delta)^{\frac{s}{2}}u_\lambda)^2\,\mathrm{d}x-\lambda\int_{\Omega}g'(x,u_\lambda)u_\lambda^2\,\mathrm{d}x\geq 0
\end{equation*}
imply that 
\begin{equation*}
\int_{\R^N}((-\Delta)^{\frac{s}{2}}u_\lambda)^2\,\mathrm{d}x=O(1)\mbox{ as }\lambda\to\Lambda^-
\end{equation*}
which yields together with monotone convergence that there exists a function $u_\Lambda$ such that $u_\lambda\to u_\Lambda$ in $\tilde{H}^s(\Omega)$ as $\lambda\to\Lambda^-$. 

Finally, we show that $(0,\Lambda)\ni\lambda\mapsto u_\lambda\in C_0(\overline{\Omega})$ is smooth. For that, we introduce the operator $A\,:\,\R^+\times C_0(\overline{\Omega})\to C_0(\overline{\Omega})$ defined by $w=A(\lambda,v)$ as the unique function satisfying
\begin{eqnarray*}
\displaystyle\left\{\begin{array}{lll}
& (-\Delta)^sw-\frac{\lambda K(x)}{w^\delta}=v\mbox{ in }\Omega,\\
& w>0\mbox{ in }\Omega,\\
& w=0\mbox{ in }\R^N\backslash\Omega.
\end{array}\right.
\end{eqnarray*}
From Theorem~\ref{sing-prob2} and Lemma~\ref{lemmaA.1}, we have
\begin{equation}\label{A-holder}
\Vert w\Vert_{C^\gamma(\R^N)}\leq C=C\left(\delta,\beta,\lambda,\Vert v\Vert_{C_0(\overline{\Omega})}\right)\quad\mbox{for some }\gamma=\gamma(\delta, \beta).
\end{equation}
From \eqref{A-holder}, we get that $A$ is a compact operator and from Lemma~\ref{lemmaA.4}, $A$ is a $C^2$ map. From Lemma~\ref{lemmaA.4} again,
the G\^ateaux derivative of $A$ at $(\lambda, v)$ in the direction $h\in C_0(\overline{\Omega})$, denoted by $\partial_2A(\lambda, v)(h)=w\in C_{\phi_{\delta,\beta}}(\Omega)$ satisfies:
\begin{equation*}
(-\Delta)^sw+\frac{\lambda\delta K(x)}{u^{\delta+1}}w=h\mbox{ with }u=A(\lambda, v).
\end{equation*}
We now define the map
\begin{equation*}
F\,:\,\R^+\times C_0(\overline{\Omega})\ni (\lambda,u)\mapsto u-A(\lambda, \lambda f(u))\in C_0(\overline{\Omega}).
\end{equation*}
From Appendix~\ref{s:Appendix_A}, $F$ is of class $C^2$ and 
\begin{equation*}
\partial_2F(\lambda,u)=I-\partial_2A(\lambda,\lambda f(u))\lambda f'(u)\mbox{ is a compact perturbation of identity}.
\end{equation*}
Since $\Lambda_1(\lambda)>0$ for any $\lambda\in (0,\Lambda)$, we obtain that $\partial_2 F(\lambda, u)$ is injective and from the Fredholm alternative, $\partial_2 F(\lambda, u)$ is invertible. Therefore, using the implicit function theorem, the map $(0,\Lambda)\ni\lambda\mapsto u_\lambda\in C_0(\overline{\Omega})$ is of class $C^2$. This completes the proof Proposition~\ref{minimal_branch}.
\qed

We now prove  the main result of this section.\\
{\bf Proof of Theorem~\ref{main-bifurcation}.}\\
We start by showing assertion i): The existence of $\Lambda$ follows from assertion i) of Proposition~\ref{minimal_branch}. Note that Theorems~\ref{sing-prob} and \ref{sing-prob2}, ${\mathcal S}\subset [0,\Lambda]\times\left (C^\gamma(\overline{\Omega})\cap C_{\phi_{\delta,\beta}}^+(\Omega)\right)$ with
\begin{eqnarray*}
\gamma=\left\{\begin{array}{lll}
& s\mbox{ if }\frac{\beta}{s}+\delta<1,\\
& s-\epsilon\mbox{ if }\frac{\beta}{s}+\delta=1\mbox{ and for }\epsilon\mbox{ small enough},\\
& \frac{2s-\beta}{\delta+1}\mbox{ if }\frac{\beta}{s}+\delta>1.
\end{array}\right.
\end{eqnarray*}
Since $2\beta+\delta(2s-1)<1+2s$, we have also that ${\mathcal S}\subset [0,\Lambda]\times \tilde{H}^s(\Omega)$. This completes the proof of assertion i). Assertion iii) follows from Proposition~\ref{minimal_branch}. Next, we prove assertion iv). If $u_\Lambda\in L^\infty$ and since $u_\lambda\uparrow u_\Lambda$ as $\lambda\uparrow\Lambda^-$, we can easily prove that $u_\Lambda\in C_{\phi_{\delta,\beta}}^+(\Omega)$. Applying the H\"older-regularity result in Theorem~\ref{sing-prob2} (see Remark~\ref{regu-problemp}), we infer that $u_\Lambda\in C^\gamma(\R^N)$ for some $\gamma=\gamma(\delta,\beta)\in (0,1)$. Therefore, $u_\lambda\uparrow u_\Lambda$ in $C_0(\overline{\Omega})$ as $\lambda\to\Lambda^-$.

Next, we apply \cite[Theorem~3.2, p.~171]{Crandall-Rabinowitz-arma} to continue the branch of minimal solutions terminating at $\Lambda$. To this aim, we note that from above
\begin{equation*}
{\rm Null}\left(\partial_2 F(\Lambda,u_\Lambda)\right)=\left\{\phi\in C_0(\overline{\Omega})\,:\, \partial_2 F(\Lambda,u_\Lambda)\phi=0\right\}\mbox{ is one dimensional and spanned by }\phi_\Lambda 
\end{equation*}
with $\phi_\Lambda>0$ and normalized in $L^2(\Omega)$. Indeed, $\Lambda_1(\Lambda)=0$. Precisely, we have proved that $\Lambda_1(\lambda)\geq 0$ for any $0<\lambda<\Lambda$. Passing to the limit as $\lambda\to\Lambda^-$, we obtain $\Lambda_1(\Lambda)\geq 0$. If $\Lambda_1(\Lambda)>0$, then the implicit function theorem asserts that the minimal branch can be continued beyond $\lambda=\Lambda$ which contradicts the definition of $\Lambda$. Therefore, $\Lambda_1(\Lambda)=0$. From the Fredholm alternative, we also have that ${\rm codim}\left({\rm Range}(\partial_2 F(\Lambda, u_\Lambda)\right)=1$. Furthermore, if $w\in {\rm Range}(\partial_2 F(\Lambda, u_\Lambda))$, then there exists $\phi\in C_0(\overline{\Omega})$ such that
\begin{equation*}
\phi-\partial_2 A(\Lambda,\Lambda f(u_\Lambda))\Lambda f'(u_\Lambda)\phi=w.
\end{equation*}
Let $w_1=\partial_2 A(\Lambda,\Lambda f(u_\Lambda))\Lambda f'(u_\Lambda)\phi$. Then, from Lemma~\ref{lemmaA.2}, $w_1\in C_{\phi_{\delta,\beta}}^+(\Omega)$ and from Lemma~\ref{lemmaA.4} satisfies
\begin{equation*}
(-\Delta)^s w_1+\frac{\Lambda\delta K(x)}{u_\Lambda^{\delta+1}}w_1=\Lambda f'(u_\Lambda)\phi.
\end{equation*}
On other hand, $\phi_\Lambda$ verifies
\begin{equation*}
(-\Delta)^s\phi_\Lambda+\frac{\Lambda\delta K(x)}{u_\Lambda^{\delta+1}}\phi_\Lambda=\Lambda f'(u_\Lambda)\phi_\Lambda.
\end{equation*}
Then, we infer that
\begin{equation}\label{w_1cara}
\int_{\Omega} wf'(u_\Lambda)\phi_\Lambda\,\mathrm{d}x=\int_{\Omega}(\phi-w_1)f'(u_\Lambda)\phi_\Lambda\,\mathrm{d}x=0.
\end{equation}
We now claim that $\partial_1 F(\Lambda, u_\Lambda)\not\in{\rm Range}(\partial_2F(\Lambda, u_\Lambda))$. Let $z=\partial_1 F(\Lambda, u_\Lambda)$. Then, from Lemma~\ref{lemmaA.2}, $z\in C_{\phi_{\delta,\beta}}^+(\Omega)$ and from Lemma~\ref{lemmaA.4} satisfies
\begin{equation*}
(-\Delta)^sz+\frac{\Lambda\delta K(x)}{u_\Lambda^{\delta+1}}z=f(u_\Lambda)+\frac{K(x)}{u_\lambda^\delta}.
\end{equation*}
Therefore, if $z\in{\rm Range}\left(\partial_2F(\Lambda, u_\Lambda)\right)$, then from \eqref{w_1cara}
\begin{equation*}
\int_{\Omega}f'(u_\Lambda)z\phi_\Lambda\,\mathrm{d}x=0
\end{equation*}
which is impossible since the integrand is strictly positive in $\Omega$. Thus \cite[Theorem~3.2, p.~171]{Crandall-Rabinowitz-arma} infers that:

Let $Z$ be a complement of ${\rm span}\{\phi_\Lambda\}$. The solutions of $F(\lambda,u)=0$ near $(\Lambda,u_\Lambda)$ are described by a curve $(\lambda(s),u(s))=(\Lambda+\tau(s), u_\Lambda+s\phi_\Lambda+x(s))$ where $s\to (\tau(s),x(s))\in\R\times Z$ is twice continuously differentiable near $s=0$ with
\begin{equation}\label{tau-0}
\tau(0)=\tau'(0)=0\mbox{ and }x(0)=x'(0)=0\mbox{ and }\tau''(0)<0.
\end{equation}
Indeed, differentiating at $s=0$, the function $F(\Lambda+\tau(s),u_\Lambda+s\phi_\Lambda+x(s))$, we get
\begin{equation*}
\tau'(0)\partial_1F(\Lambda, u_\Lambda)+\partial_2F(\Lambda,u_\Lambda)(\phi_\Lambda+x'(0))=0.
\end{equation*}
Noticing that 
\begin{equation*}
\partial_2F(\Lambda,u_\Lambda)\phi_\Lambda=0\,\mbox{ and} \,\partial_1F(\Lambda,u_\Lambda)\not\in\left({\rm Range}(\partial_2 F(\Lambda, u_\Lambda)\right),
\end{equation*}
we get $x'(0)=0$ and $\tau'(0)=0$. 
Differentiating again the function $F(\lambda(s), u(s))$, with respect to $s$ and evaluating at $s=0$, we obtain:
\begin{equation}\label{eq*}
\tau''(0)\partial_1 F(\Lambda, u_\Lambda)+\partial_2F(\Lambda,u_\Lambda)x''(0)+\partial^2_{22}F(\Lambda, u_\Lambda)(\phi_\Lambda,\phi_\Lambda)=0.
\end{equation}
Thus, 
\begin{equation*}
\tau''(0)z+\partial_2F(\Lambda,u_\Lambda)x''(0)+\tilde{w}=0
\end{equation*}
where $\tilde{w}=\partial^2_{22}F(\Lambda,u_\Lambda)(\phi_\Lambda,\phi_\Lambda)$ and by Lemma~\ref{lemmaA.4} solves the equation:
\begin{equation*}
(-\Delta)^s\tilde{w}+\frac{\Lambda\delta K(x)}{u_\Lambda^{\delta+1}}\tilde{w}=\frac{\lambda\delta(\delta+1)K(x)}{u_\Lambda^{\delta+2}}\phi_\Lambda^2+\Lambda f''(u_\Lambda)\phi_\Lambda^2>0\mbox{ in }\Omega
\end{equation*}
since by {\bf (f3)}, $t\to g(x,t)=K(x)t^{-\delta}+f(t)$ is strictly convex.

Multiplying the equation in \eqref{eq*} by $\phi_\Lambda$ and integrating on $\Omega$, we obtain:
\begin{equation}\label{tauseconde}
\tau''(0)\int_{\Omega}z\phi_\Lambda\,\mathrm{d}x+\int_{\Omega}\partial_2F(\Lambda,u_\Lambda)x''(0)\phi_\Lambda\,\mathrm{d}x=-\int_{\Omega}\tilde{w}\phi_\Lambda\,\mathrm{d}x.
\end{equation}
Recalling that the middle term in \eqref{tauseconde} vanishes, we obtain that $\tau''(0)<0$. This completes the proof of assertion iv). let us finally prove the assertion ii).

The existence of a connected unbounded branch of solutions to $(P_\lambda)$ can be proved similarly as in the proof of \cite[Theorem~3.2, p.~508]{Ra-JFA} by a Leray-Schauder argument.

Assertion v) follows from the fact that there is no solution to $(P_\lambda)$ for $\lambda>\Lambda$ together with the unboundedness of ${\mathcal C}$. This completes the proof of Theorem~\ref{main-bifurcation}.
\qed
\section{Applications}\label{section3}
In this section, we prove Theorem~\ref{bifurcation-example} and Theorem~\ref{example-analytic}. We start with\\
{\bf Proof of Theorem~\ref{bifurcation-example}.}\\
We first prove that $u_\Lambda\in L^\infty$. For that we use the subcritical growth of $f$ (from ${\bf (f4)}$) and the following inequality :
\begin{eqnarray*}
(-\Delta)^s(u_\Lambda-1)^+\leq\Lambda( C+f(u_\Lambda))\leq C_0\Lambda(1+u_\Lambda^p)\mbox{ in }\Omega
\end{eqnarray*}
which holds for some constants $C, C_0>0$. Then, we can apply a classical bootstrap arguments together with regularity results in \cite[Proposition 1.4, p.~727]{Ros-othon-serra-CVPDE} to get $u_\Lambda\in L^\infty$.

From Theorem~\ref{main-bifurcation}, it is now sufficient to prove uniform estimates in ${\mathcal S}\cap\left(\{\lambda\geq \lambda_0\}\times C_0(\overline{\Omega})\right)$ for any $\lambda_0>0$.

For that, we use a similar approach as in the proof of \cite[theorem~1, p.~148]{Chen-Li-Arma} (see also a priori estimates in the same fashion in \cite{Figueiredo-Lions_Nussbaum-JMPA}). Precisely, we prove uniform estimates near the boundary by using the moving plane method whereas interior estimates are derived through a blow up analysis together with a suitable Liouville theorem as in \cite{Gidas-Spruck}.

In this regard, from\cite{Chen-Li-Arxiv}, we have the following Liouville type result:\\
Let $u\in C^{1,1}_{\rm loc}(\R^N)$, nonnegative satisfying $\int_{\R^N}\frac{u}{1+\vert x\vert^{N+2s}}\,\mathrm{d}x<\infty$ and
\begin{equation*} 
(-\Delta)^s u=u^p\;\mbox{ in }\R^N, \;\mbox{ with }1<p<\frac{N+2s}{N-2s}.
\end{equation*}
Then, $u\equiv 0$.

To perform the moving plane method, we need a main ingredient : a maximum principle for narrow domains. We argue as in Section 2.2 in \cite{Chen-Li-Arma}. Let $x_0\in \partial\Omega$. We suppose first that $\Omega$ is strictly convex in a neighborhood of $x_0$. Without loss of generality, we can assume that the outward normal at $x_0$, $\boldsymbol{\nu}(x_0)$, satisfies $\boldsymbol{\nu}(x_0)=(-1,0,\cdot\cdot,0)$. We then define
\begin{equation*}
T_\mu=\left\{x\in\R^N\,:\, x_1=\mu\right\}\mbox{ for some }\mu\in\R,
\end{equation*}
\begin{equation*}
\Sigma_\mu=\left\{x\in\R^N\,:\, x_1<\mu\right\}
\end{equation*}
and $x^\mu=(2\mu-x_1,x_2,\cdot,\cdot, x_N)$ be the reflection point of $x=(x_1,x_2,\cdot,\cdot, x_N)$ about the plane $T_\mu$.

Let $\lambda_0>0$ and let $u$ be  a solution to $(P_\lambda)$ with $\lambda\geq\lambda_0$. we compare the values of $u(x)$ with $u^\mu(x)\eqdef u(x^\mu)$. For that, we denote $w^\mu(x)\eqdef u^\mu(x)-u(x)$. For $\mu$ sufficiently negative, we have clearly $w^\mu(x)\geq 0$ in $\Sigma_\mu$. We need to prove $w^\mu\geq 0$ when $T_\mu$ meets a neighborhood of $x$. In thsi step,  we use a maximum principle for narrow domains: 
Suppose that $u^\mu<0$ in a region $D\subset \Sigma_\mu$. Then,
\begin{equation*}
\left\langle(-\Delta)^s(-w^\mu),(-w^\mu)^+\right\rangle=\lambda\int_{\Omega}\left(\frac{1}{u^\delta}-\frac{1}{(u^\mu)^\delta}+f(u)-f(u^\mu)\right)(-w^\mu)^+\,\mathrm{d}x.
\end{equation*}
Then, using that $f$ is a Lipschitz function, we obtain for a constant $C>0$ :
\begin{equation*}
\int_{\R^N}\left((-\Delta)^{\frac{s}{2}}(u-u^\mu)^+\right)^2\,\mathrm{d}x\leq C\int_{D}((u-u^\mu)^+)^2\,\mathrm{d}x
\end{equation*}
and by the Poincar\'e inequality (with the associated best constant $C_p$):
\begin{equation*}
\int_{\R^N}\left((-\Delta)^{\frac{s}{2}}(u-u^\mu)^+\right)^2\,\mathrm{d}x\leq CC_p(D)\int_{\R^N}\left((-\Delta)^{\frac{s}{2}}(u-u^\mu)^+\right)^2\,\mathrm{d}x.
\end{equation*}
If the diameter of $D$ is small enough, then $CC_p(D)<1$ and then $(u-u^\mu)^+\equiv 0$. Alternatively, we can use the maximum principle for narrow domains  \cite[Theorem~2, p.~7]{Chen-Li-Arxiv}. We apply Theorem~2 by setting the function $c(x)=\frac{\lambda\delta}{(u+\theta (u^\mu-u))^{\delta+1}}(x)+\frac{f(u)-f(u^\mu)}{u^\mu-u}(x)$, with $\theta\in (0,1)$, is bounded by below since $f$ is Lipschitz. Note that the quantity referred as $\delta$ in  \cite[Theorem~2, p.~7]{Chen-Li-Arxiv} is not depending of $w^\mu=u^\mu-u$ (see proof of  \cite[Theorem~2.3, p.~13]{Chen-Li-Arxiv}). 

Now, we are in the following situation: By moving the hyperplane in a direction close to the outward normal in a neighborhood of any $x_0\in\partial\Omega$, we infer that there exist a $H>0$ and a $T>0$ independent of $u$ such that $u(x-t\gamma)$ is non increasing for $t\in [0,T]$ and for any $x$ in a neighborhood of $\tilde{x}$ and for any $\gamma\in \R^N$ satisfying $|\gamma |=1$ and $\gamma\cdot\boldsymbol{\nu}(\tilde{x})\geq H$ for all $\tilde{x}\in\partial\Omega$.

The fact that $u(x-t\gamma)$ is non increasing in $t$ for $x$ and $\gamma$ described above implies that we have two positive numbers $\alpha_1$ and $\alpha_2$ both depending on $\partial\Omega$ such that for any $x\in\Omega_{\alpha_2}\eqdef\{x\in\Omega\,:\,d(x)<\alpha_2\}$, we have a measurable set $I_x$ verifying:
\begin{itemize}
\item[(i)] $\vert I_x\vert\geq \alpha_1$;
\item[(ii)] $I_x\subset\{x\in\Omega\,:\, d(x)\geq\frac{\alpha_2}{2}\}$;
\item[(iii)] $u(y)\geq u(x),$ $\forall y\in I_x$.
\end{itemize}
From above, we deduce a uniform a priori bound in a neighborhood of $\partial\Omega$: Multiplying by $\phi_{1,s}$ the equation satisfied by $u$ we get:
\begin{equation}\label{eq2.2}
\lambda_{1,s}\int_{\Omega} u\phi_{1,s}\,\mathrm{d}x=\lambda\left(\int_{\Omega}\frac{\phi_{1,s}}{u^\delta}\,\mathrm{d}x+\int_{\Omega}f(u)\phi_{1,s}\,\mathrm{d}x\right).
\end{equation}
Observing from {\bf (f3)} that for any $\rho>\frac{\lambda_{1,s}}{\lambda_0}$, there exists $C>0$ such that
\begin{equation}\label{eq2.3}
\frac{1}{t^\delta}+f(t)\geq \rho t-C\quad\forall t\in\R^+
\end{equation}
and using \eqref{eq2.2}, it follows that for some constant $C>0$
\begin{equation*}
(\rho-\frac{\lambda_{1,s}}{\lambda_0})\int_{\Omega}u\phi_{1,s}\,\mathrm{d}x\leq C
\end{equation*}
which implies for $\tilde{C}\eqdef \frac{C}{\rho-\frac{\lambda_{1,s}}{\lambda_0}}$
\begin{equation*}
u(x)\int_{I_x}\phi_{1,s}\,\mathrm{d}x\leq\int_\Omega u\phi_{1,s}\,\mathrm{d}x\leq \tilde{C}.
\end{equation*}
Thus, since $\phi_1(x)\geq c_2d(x)$ for some $c_2>0$ and taking into account (i) and (iii)
\begin{equation*}
u(x)\leq \frac{\tilde C}{\alpha_1\alpha_2c_2}\mbox{ for all }x\in \Omega_{\alpha_2}.
\end{equation*}
This completes the proof of uniform estimates near the boundary. 

Next we prove the uniform interior estimates by a blow-up technique as in \cite{Gidas-Spruck}. Precisely, suppose that there exists a sequence $(\lambda_k,u_k)_{k\in\N^*}$ of solutions in ${\mathcal S}$ such that for some $\lambda_0>0$, $\lambda_0\leq \lambda_k\leq \Lambda$ and $\Vert u_k\Vert_{L^\infty(\Omega)}\to\infty$ as $k\to\infty$. Let $x_k={\rm argmax}(u_k)$ and $M_k=u_k(x_k)=\displaystyle\max_{\Omega}u_k$. Note that from the uniform estimates near the boundary established above, there exists $c_0>0$ such that
\begin{equation*}
d(x_k)\geq c_0\quad\forall k\in \N^*.
\end{equation*}
Then, we define the rescaled function $v_k$ for all $k\in\N^*$ such as:
\begin{equation*}
v_k(y)=\mu_k^{\frac{2s}{p-1}}u_k(x),\quad y=\frac{x-x_k}{\mu_k}\in\Omega_k\eqdef\frac{\Omega-x_k}{\mu_k}\;\mbox{ and }\;\mu_k^{\frac{2s}{p-1}}M_k=1.
\end{equation*}
Consequently, $v_k$ satisfies
\begin{equation*}
(-\Delta)^sv_k(y)=\lambda_k\left(\frac{\mu_k^{\frac{2(p+\delta)s}{p-1}}}{v_k^\delta(x)}+\mu_k^{\frac{2ps}{p-1}}f(\mu_k^{\frac{-2s}{p-1}}v_k)\right)\mbox{ in }\Omega_k.
\end{equation*}
As in \cite{Gidas-Spruck}, we can prove that up to a subsequence, $v_k\to v$ in $C^s_{\rm loc}(\R^N)$ as $k\to\infty$ and from {\bf (f4)} $v$ satisfies:
\begin{eqnarray*}
\left\{\begin{array}{ll}
&(-\Delta)^s v= cv^p\mbox{ in }\R^N,\\
&v(0)=1.
\end{array}\right.
\end{eqnarray*}
From assertion (ii) of \cite[Theorem~4, p.~8]{Chen-Li-Arxiv} (Liouville theorem), we get a contradiction. Therefore we obtain a uniform $L^\infty$-bound of solutions in ${\mathcal S}\cap\{\lambda\geq\lambda_0\}$, for any $\lambda_0>0$. From assertion (v) of Theorem~\ref{main-bifurcation}, we then conclude that $0$ and only $0$ is an asymptotic bifurcation point.

Finally, we deal with the general case, i.e. where $\Omega$ is not  strictly convex. In this case as in \cite{Chen-Li-Arma}, we perform a Kelvin transform near any boundary point $x_0$. Precisely, let $K_0=\displaystyle\max_{x\in\partial\Omega} k(x)$ where $k(x)$ denotes the curvature of $\partial\Omega$ at $x\in\partial\Omega$. Consider $R=K_0+1$, $x_1=x_0+\frac{\boldsymbol{\nu}(x_0)}{R}$ and assuming that the outward normal of $\Omega$ at $x_0$ $\boldsymbol{\nu}(x_0)=(-1,0,\cdot,\cdot,0)$. Then, $B_{1/R}(x_1)$ is tangent to $\Omega$ at $x_0$ and $B_{1/R}(x_1)\cap \Omega=\emptyset$. We use the following inversion transformation $T$:
\begin{equation*}
T:\, x\mapsto Tx\eqdef\frac{x-x_1}{|x-x_1|^2}=y\mbox{ and }x=x_1+\frac{y}{|y|^2}.
\end{equation*}
We have that $T(\Omega)\subset B_R(0)$ and $T(\Omega)$ is tangent to $B_R(0)$ at $(R,0,\cdot\cdot,0)$. Next, we define $u^*$ by 
\begin{equation*}
u^*(y)=|y|^{2s-N}u(x_1+\frac{y}{|y|^2}).
\end{equation*}
Then, we have
\begin{eqnarray*}
(-\Delta)^su^*(y)&=|y|^{-2s-N}(-\Delta)^su(x)=|y|^{-2s-N}\left(\frac{\lambda}{u^\delta(x)}+\lambda f(u(x))\right)\\
&=\frac{\lambda}{|y|^{2s+N+(N-2s)\delta}u^*(y)^\delta}+\frac{\lambda f( |y|^{N-2s} u^*(y))}{|y|^{2s+N}}.
\end{eqnarray*}
As in subsection 2.2 in \cite{Chen-Li-Arma}, we need to prove that
\begin{equation*}
\frac{\partial}{\partial y}f^*(y,u^*)\leq 0 \mbox{ where }f^*(y,u^*)=\frac{\lambda}{|y|^{2s+N+(N-2s)\delta}(u^*)^\delta}+\frac{\lambda f( |y|^{N-2s} u^*)}{|y|^{2s+N}}.
\end{equation*}
It is sufficient to verify that
\begin{equation*}
\frac{\partial}{\partial y}\left(\frac{\lambda f( |y|^{N-2s} u^*(y))}{|y|^{2s+N}}\right)\leq 0.
\end{equation*}
From {\bf (f5)}, we have
\begin{eqnarray*}
&\frac{\partial}{\partial y}\left(\frac{\lambda f( |y|^{N-2s} u^*(y))}{|y|^{2s+N}}\right)=\frac{-(N+2s)}{|y|^{N+2s+1}}f([y|^{N-2s}u^*)+\frac{(N-2s)|y|^{N-2s-1}u^*}{|y|^{N+2s}}f'(|y|^{N-2s}u^*)\\
&\leq\frac{f(|y|^{N-2s}u^*)}{|y|^{N+2s+1}}\left(-N-2s+q(N-2s)\right)\leq 0.
\end{eqnarray*}
This completes the proof of Theorem~\ref{bifurcation-example}.
\qed
\begin{rema}\label{global-mult}
From Theorem~\ref{bifurcation-example}, we get a global multiplicity result that extends results in \cite{peral-al}.
\end{rema}
\begin{rema}\label{caseN=1}
Similar results as in Theorem~\ref{bifurcation-example} can be derived in the case $N=1$, $s=\frac{1}{2}$, and considering $f$ with at most exponential growth. The boundary estimate can be proved via the moving plane method whereas the interior estimates can be performed using the blow-up analysis of Brezis-Merle type proved in \cite[Theorem~1.1, p.~1758]{Da Lio-Martinazzi-Riviere}, the Trudinger-Moser inequality in \cite[Theorem~1, p.~264]{Martinazzi} (see also \cite{GiPaSe}) and the expression and behaviour of the fundamental solution (see \cite[Theorem~3.1, p.~26]{Bucur} or \cite[Corollary 1.2, p.~1309]{ChSoKi}). We will discuss this case in details in a forthcoming paper.
\end{rema}

Now, we consider the particular case $f(u)=u^p$ with $1<p<\frac{N+2s}{N-2s}$ and $\frac{\beta}{s}+\delta<1$. We prove in this case the existence of an analytic branch of solutions to $(P_\lambda)$. For that we appeal the analytic global bifurcation theory introduced by Dancer (see \cite{Dancer}). Precisely we apply a variant form of \cite[Theorem~9.1.1, p.~114]{Buffoni-Toland-book}, see also \cite{GiPrWa}). To start with, recall that the operator $F\,:\,\R^+\times C_{\phi_{1,s}}^+(\Omega)\mapsto C_{\phi_{1,s}}^+(\Omega)$ is defined by
\begin{equation*}
F(\lambda, u)\eqdef(-\Delta)^{-s}\left(\frac{\lambda K(x)}{u^\delta}+\lambda u^p\right)\mbox{ for any }(\lambda,u)\in \R^+\times C_{\phi_{1,s}}^+(\Omega).
\end{equation*}
To prove the $0$-index Fredholmness of the operator $I-\partial_2 F$, we establish the following compactness lemma. The lemma below also ensures the compactness of bounded closed subsets of ${\mathcal S}$ in $\R^+\times C_{\phi_{1,s}}(\Omega)$ and implies Lemma~\ref{Compactness-Lemma}.
\begin{lemm}\label{compactness-lemma}
Let $u\in C_{\phi_{1,s}}^+(\Omega)$ and consider the operator $\tilde{T}\,:\, C_{\phi_{1,s}}(\Omega)\mapsto C_{\phi_{1,s}}(\Omega)$ defined by
\begin{equation*}
\tilde{T}(\phi)\eqdef (-\Delta)^{-s}\left(\frac{K(x)\phi}{u^{\delta +1}}\right),\quad\forall\phi\in C_{\phi_{1,s}}(\Omega).
\end{equation*}
Then $\tilde{T}$ is compact.
\end{lemm}
{\bf Proof of Lemma~\ref{compactness-lemma}.}\\
Let $(w_n)_{n\in\N}\subset C_{\phi_{1,s}}(\Omega)$ be a bounded sequence in $C_{\phi_{1,s}}(\Omega)$, i.e. satisfying $\displaystyle\sup_{n\in\N^*}\left\Vert \frac{w_n}{d^s}\right\Vert_{L^\infty(\Omega)}<\infty$. Let $v_n=(-\Delta)^{-s}\left(\frac{w_n K(x)}{u^{\delta+1}}\right)$. We will show that $(v_n)_{n\in\N}$ is relatively compact in $C_{\phi_{1,s}}(\Omega)$. For that, let $\epsilon>0$ and define $Z_\eta=\left\{x\in \Omega\,:\, d(x)\geq \eta\right\}$, the corresponding indicator function $\un_{Z_\eta}$ and
\begin{eqnarray*}
&v_n^{1,\epsilon}\eqdef(-\Delta)^{-s}\left(\frac{w_n K(x)\un_{Z_\epsilon}}{u^{\delta +1}}\right),\quad v_n^{2,\epsilon}\eqdef\left((-\Delta)^{-s}\left(\frac{w_n K(x)(1-\un_{Z_{\epsilon}})}{u^{\delta +1}}\right)\right)\un_{Z_{3\epsilon}}\mbox{ and }\\
&v_n^{3,\epsilon}\eqdef\left((-\Delta)^{-s}\left(\frac{w_n K(x)(1-\un_{Z_{3\epsilon}})}{u^{\delta +1}}\right)\right)(1-\un_{Z_{3\epsilon}}).
\end{eqnarray*}
Clearly, $v_n=v_n^{1,\epsilon}+v_n^{2,\epsilon}+v_n^{3,\epsilon}$. So it is sufficient to prove that for $i=1,2,3$, $(v_n^{i,\epsilon})_{n\in\N}$ is relatively compact in $C_{\phi_{1,s}}(\Omega)$.

We first prove that for $\epsilon>0$ fixed, $(v_n^{1,\epsilon})_{n\in\N}$ and $(v_n^{2,\epsilon})_{n\in\N}$ are relatively compact in $C_{\phi_{1,s}}(\Omega)$. Concerning the sequence $(v_n^{1,\epsilon})_{n\in\N}$, we observe that $\frac{w_n K(x)\un_{Z_\epsilon}}{u^{\delta +1}}\in L^\infty(\Omega)$ and from \cite[Theorem~1.2,p.~277]{Ros-oton-serra-JMPA}, we obtain that 
\begin{equation*}
\left\Vert\frac{v_n^{1,\epsilon}}{d^s}\right\Vert_{C^s(\R^N)}\leq C=C(\epsilon).
\end{equation*}
Therefore, for $\epsilon>0$ fixed, $(v_n^{1,\epsilon})_{n\in\N}$ is relatively compact in $C_{\phi_{1,s}}(\Omega)$. let us now consider the sequence $(v_n^{2,\epsilon})_{n\in\N}$. For any $x,x'\in Z_{3\epsilon}$,
\begin{eqnarray*}
\left\vert \frac{v_n^{2,\epsilon}(x)}{d^s(x)}- \frac{v_n^{2,\epsilon}(x)}{d^s(x')}\right\vert &=\left\vert\int_{\Omega}\left(\frac{G(x,y)}{d^s(x)}-\frac{G(x',y)}{d^s(x')}\right)\frac{K(y)(1-\un_{Z_\epsilon})w_n(y)}{u^{\delta +1}(y)}\,\mathrm{d}y\right\vert\\
&\leq C\int_{\Omega}\left\vert \frac{G(x,y)}{d^s(x)}-\frac{G(x',y)}{d^s(x')}\right\vert\frac{(1-\un_{Z_\epsilon}(y))}{d^{\beta+\delta s}(y)}\,\mathrm{d}y.
\end{eqnarray*}
Where $G(x,y)$ denotes the Green's function associated to $(-\Delta)^s$ in $\Omega$ with Dirichlet boundary conditions.
Next, we prove that $x\mapsto\frac{G(x,y)}{d^s(x)}$ is H\"older-continuous in $Z_{3\epsilon}$ uniformly with respect to $y\in \Omega\backslash Z_{\epsilon}$ (but still depending on $\epsilon$). Using the estimate in Corollary~\ref{coro-holder} and a finite balls covering, we deduce that
\begin{equation*}
\Vert G(x,y)\Vert_{C^s(K_{3\epsilon})}\leq C\left(\Vert (1+\vert x\vert)^{-N-2s} G(x,y)\Vert_{L^1(\R^N)}+\Vert G(x,y)\Vert_{L^\infty(K_{2\epsilon})}\right).
\end{equation*}
Furthermore, for any $y\in\R^N\backslash Z_\epsilon$ and a fixed $R=R(\epsilon)>0$ small enough, there exists $C=C(\epsilon)>0$ such that
\begin{equation*}
\Vert G(x,y)\Vert_{L^\infty(Z_{2\epsilon})}\leq C\mbox{ for any }y\in \R^N\backslash Z_\epsilon\mbox{ and }
\end{equation*}
\begin{equation*}
\int_{\R^N}\frac{1}{(1+\vert x\vert)^{N+2s}}G(x,y)\,\mathrm{d}x\leq C\int_{B_R(y)}\frac{\mathrm{d}x}{\vert x-y\vert^{N-2s}}+C\int_{\R^N\backslash B_R(y)}\frac{\mathrm{d}x}{(1+\vert x\vert)^{N+2s}}<\infty.
\end{equation*}
Therefore, there exists $C_\epsilon>0$ such that $\Vert G(x,y)\Vert_{C^s(Z_{3\epsilon})}\leq C_\epsilon$ uniformly with respect to $y\in \R^N\backslash Z_\epsilon$. Then, we deduce that for  some constant $\tilde{C}_\epsilon$,
\begin{equation*}
\left\vert\frac{v_n^{2,\epsilon}(x)}{d^s(x)}-\frac{v_n^{2,\epsilon}(x')}{d^s(x')}\right\vert\leq C_\epsilon\vert x-x'\vert^s\int_{\Omega}\frac{\mathrm{d}y}{d(y)^{\beta+\delta s}}\leq \tilde{C}_\epsilon\vert x-x'\vert^s
\end{equation*}
from which it follows that for a fixed $\epsilon>$, $(v_n^{2,\epsilon})_{n\in\N}$ is relatively compact in $C_{\phi_{1,s}}(\Omega)$.

Now, we will prove uniform estimates depending on $\epsilon$ on the sequence $(v_n^{3,\epsilon})_{n\in\N}$. Let $\gamma'$ such that $\beta+ s\delta<\gamma'<s$. Then,
\begin{eqnarray*}
&\left\vert\frac{(1-\un_{Z_{3\epsilon}})}{d^s(x)}\int_{\R^N\backslash Z_\epsilon}\frac{G(x,y)}{d^{\beta+s\delta}(y)}\,\mathrm{d}y\right\vert\leq  C\frac{(1-\un_{Z_{3\epsilon}})}{d^s(x)}\int_{\R^N\backslash Z_\epsilon}\displaystyle\min\left(\frac{d^s(x)d^s(y)}{\vert x-y\vert^N},\frac{d^s(x)}{\vert x-y\vert^{N-s}}\right)\frac{\mathrm{d}y}{d(y)^{\beta+\delta s}}\\
&\leq  C (1-\un_{Z_{3\epsilon}})\epsilon^{\gamma'-(\beta+\delta s)}\int_{\Omega}\displaystyle\min\left(\frac{d(y)^s}{\vert x-y\vert^N},\frac{1}{\vert x-y\vert^{N-2s}}\right)\frac{\mathrm{d}y}{d(y)^{\gamma'}}\leq O(\epsilon^{\gamma'-(\beta+s\delta)}).
\end{eqnarray*}
We finally prove that $\left(\frac{v_n}{d^s}\right)_{n\in\N}$ is relatively compact in $L^\infty(\Omega)$. Let $\eta>0$ be small enough. First fix $\epsilon>0$ small enough such that $\displaystyle\left\Vert\frac{v_n^{3,\epsilon}}{d^s}\right\Vert_{L^\infty(\Omega)}\leq \eta$. Then, for such $\epsilon$, we can extract convergent subsequences $(v_{\psi(n)}^{1,\epsilon})_{n\in\N}$ and $(v_{\psi(n)}^{2,\epsilon})_{n\in\N}$ in $C_{\phi_{1,s}}(\Omega)$. Therefore, there exists $M=M(\eta)\in \N$ large enough such that for $n,m\geq M$
\begin{eqnarray*}
\left\Vert \frac{v_{\psi(n)}}{d^s}-\frac{v_{\psi(m)}}{d^s}\right\Vert_{L^\infty(\Omega)}\leq\left\Vert \frac{v_{\psi(n)}^{1,\epsilon}}{d^s}-\frac{v_{\psi(m)}^{1,\epsilon}}{d^s}\right\Vert_{L^\infty(\Omega)}+\left\Vert \frac{v_{\psi(n)}^{2,\epsilon}}{d^s}-\frac{v_{\psi(m)}^{2,\epsilon}}{d^s}\right\Vert_{L^\infty(\Omega)}+\eta\leq3\eta.
\end{eqnarray*}
Therefore, $(v_{\psi(n)})_{n\in\N}$ is a Cauchy sequence in $C_{\phi_{1,s}}(\Omega)$ and then $(v_n)_{n\in\N}$ is relatively compact in $C_{\phi_{1,s}}(\Omega)$.
\qed

From Lemma~\ref{compactness-lemma}, we deduce the following corollary:
\begin{coro}
$A\,:\,\R^+\times C_{\phi_{1,s}}^+(\Omega)\mapsto C_{\phi_{1,s}}^+(\Omega)$ and $\partial_2 A(\lambda,u)$, for $(\lambda,u)\in \R^+\times C_{\phi_{1,s}}^+(\Omega)$  are compact operators.
\end{coro}

We are now ready to give\\
{\bf Proof of Theorem~\ref{example-analytic}.}\\
Using similar arguments as in \cite[Proposition~1,p.~372]{Da2}, we can prove that
\begin{equation*}
F\,:\,\R^+\times C_{\phi_{1,s}}^+(\Omega)\mapsto C_{\phi_{1,s}}^+(\Omega)
\end{equation*}
is real analytic. In particular, we can prove that
\begin{equation*}
G\,:\, C_{\phi_{1,s}}^+(\Omega)\mapsto C_{\phi_{1,s}^{-\frac{\beta}{s}-\delta}}^+(\Omega)
\end{equation*}
defined by $G(u)=K(x)u^{-\delta}+u^p$ is analytic and
\begin{equation*}
(-\Delta)^{-s}\,:\, C_{\phi_{1,s}^{-\frac{\beta}{s}-\delta}}(\Omega)\mapsto C_{\phi_{1,s}}(\Omega)
\end{equation*}
is a linear continuous map and maps  $C_{\phi_{1,s}^{-\frac{\beta}{s}-\delta}}^+(\Omega)$ into $C_{\phi_{1,s}}^+(\Omega)$.  Defining
the non-singular solution set
$$\mathcal{N}= \displaystyle\left\{(\lambda,x)\in \mathcal{S}\,:\,{\rm Null}( \partial_2 F(\lambda,x))=\{0\}\right\}$$
and a distinguished arc as a maximal connected subset of ${\mathcal N}$, we  state below the global bifurcation result in the analytic framework that we use as  a variant of \cite[Theorem~9.1.1]{Buffoni-Toland-book}:

Suppose that
\begin{itemize}
\item[(G1)] Bounded closed subsets of ${\mathcal S}$ are compact in $\R \times \mathcal{X}$.
\item[(G2)]  $\partial_2F(\lambda, x)$ is a Fredholm operator of index zero for all $(\lambda,x) \in \mathcal{S}$. 
\item[(G3)] There exists an analytic function $(\lambda,u) \, : (0,\epsilon) \to \mathcal{S}$ such that     $\partial_2F(\lambda(s),u(s))$ is invertible for all $s\in (0,\epsilon)$ and  $\displaystyle\lim_{s\to 0^+}(\lambda(s), u(s))=(0,0)$.
\end{itemize}
 Let
$$ {\mathcal A_0}=\left\{(\lambda(s),u(s))\,:\, s\in (0,\epsilon)\right\}.$$
Obviously, ${\mathcal A_0}\subset {\mathcal S}$. The following result gives a global extension of the function $(\lambda, u)$ from $(-\epsilon,\epsilon)$ to $(-\infty,\infty)$ in the real analytic case.
\begin{theo}\label{theo9.1.1}
Suppose (G1)-(G3) hold. Then,  $(\lambda,u)$ can be extended as a continuous map (still called) $(\lambda,u) :  (-\infty,\infty) \to  \mathcal{S}$ with the following properties:
\begin{itemize}
\item[(a)] Let ${\mathcal A} \eqdef \{(\lambda(s),u(s)): s \in\R\}.$ Then, ${\mathcal A} \cap {\mathcal N}$ is an atmost countable union of distinct distinguished arcs $\bigcup_{i=1}^n {\mathcal A_i},\; n \leq \infty$.
\item[(b)] ${\mathcal A_0}\subset {\mathcal A_1}$.
\item[(c)]  $\{s\in\R\,:\,ker(\partial_2 F(\lambda(s),u(s))) \neq \{0\}\}$ is a discrete set.
\item[(d)] At each of its points ${\mathcal A}$ has a local analytic re-parameterization in the following sense: For each $s^*\in \R$ there exists a continuous, injective map $\rho^*\,:\, (-1,1)\to \R$ such that $\rho^*(0)=s^*$  and the re-parametrisation
\begin{eqnarray*}
 (-1,1) \ni t\to (\lambda(\rho^*(t)),u(\rho^*(t))) \in \mathcal{A} \mbox{ is analytic}.
\end{eqnarray*}
Furthermore, the map $s \mapsto \lambda(s)$ is injective in a  neighborhood of $s=0$ and for each $s^*\neq 0$ there exists $\epsilon^*>0$ such that $\lambda$ is injective on $[s^*,s^*+\epsilon^*]$ and on $[s^*-\epsilon^*,s^*]$.
\item[(e)] Only one of the following alternatives occurs:
\begin{itemize}
\item[(i)] $\Vert(\lambda(s),u(s))\Vert_{\R \times \mathcal{X}}\to\infty$ as $s\to +\infty$ (resp. $s \to -\infty$).
\item[(ii)] a subsequence $\{(\lambda(s_n),u(s_n))\}$ approaches the boundary of $\mathcal{U}$ as $s_n \to +\infty$ (resp. $s_n \to -\infty$).
\item[(iii)] ${\mathcal A}$ is the closed loop :
$$
{\mathcal A}=\left\{(\lambda(s),u(s))\,:\, -T\leq s\leq T, (\lambda(T),u(T))=(\lambda(-T),u(-T)) \text{ for some } T>0 \right\}.
$$ In this case, choosing the smallest such $T>0$ we have
\begin{eqnarray*}
(\lambda(s+2T),u(s+2T))=(\lambda(s),u(s)) \mbox{ for all } s \in \R.
\end{eqnarray*}
\end{itemize}
\item[(f)] Suppose $\partial_2 F(\lambda(s_1),u(s_1))$  is invertible  for some $s_1\in\R$. If for some $s_2\neq s_1$, we have 
$(\lambda(s_1),u(s_1))=(\lambda(s_2),u(s_2))$
then (e)(iii) occcurs and $\vert s_1-s_2\vert$ is an integer multiple of $2T$. In particular, the map $s \mapsto (\lambda(s), u(s))$ is  injective on $[-T,T)$.
\end{itemize}
\end{theo}
Let us check that the assumptions below are satisfied. Conditions (G1), (G2) follow from Lemma~\ref{compactness-lemma} with ${\mathcal X}=C_{\phi_{1,s}}(\Omega)$. From the analytic version of the implicit function theorem and Proposition~\ref{minimal_branch},
\begin{equation*}
{\mathcal A_0}=\{(\lambda,u_\lambda)\,:\,\lambda\in (0,\Lambda)\}
\end{equation*}
satisfies statements in (G3). We fix an analytic parametrization ${\mathcal A_0}=\{(\lambda(s),u(s))\,:\, s\in (0,s_0)\}$ for some $s_0>0$.  Applying the version of 
\cite[Theorem~9.1.1]{Buffoni-Toland-book} above, we infer the extension of the  analytic map $(\lambda(s),u(s)$ for all $s>0$ to get a global analytic and continuous branch ${\mathcal A}$ of solutions to $(P_\lambda)$ containing ${\mathcal A_0}$ and satisfying assertions (a)-(f). This proves i), ii), v) and vi) of Theorem~\ref{example-analytic}. The assertion e(i) occurs since the branch ${\mathcal A_0}$ emanating from $(0,0)$ is unique in its some neighborhood and since ${\mathcal A}\subset [0,\Lambda]\times C_{\phi_{1,s}}^+(\Omega)$. Then, assertion iv) follows from the unboundedness of ${\mathcal A}$ and assertion vii) is a consequence of the nonexistence of bounded weak solutions for $\lambda>\Lambda$. This completes the proof of Theorem~\ref{example-analytic}.
\qed

\appendix
\section{}
\label{s:Appendix_A}
In this appendix, assuming that $2\beta+\delta(2s-1)<1+2s$, we prove the regularity of $A\,:\, \R^+\times C_0(\overline{\Omega})\ni(\lambda,h)\mapsto u\in C_0(\overline{\Omega})$, defined as
$A(\lambda,h)=u$ unique solution to
\begin{eqnarray}\label{eqA1}
\left\{\begin{array}{lll}
&(-\Delta)^s u-\frac{\lambda K(x)}{u^\delta}=h\mbox{ in }\Omega\\
& u>0\mbox{ in }\Omega\\
& u=0\mbox{ in }\R^N\backslash\Omega
\end{array}\right.
\end{eqnarray}
About regularity of $v$, we have the following lemma:
\begin{lemm}\label{lemmaA.1}
Let $h\in C_0(\overline{\Omega})$ and $\lambda\in\R^+$. Then, $u=A(\lambda, h)\in \tilde{H}^s(\Omega)\cap C_{\phi_{\delta,\beta}}^+(\Omega)\cap C^\gamma(\R^N)$ with $\gamma=\gamma(\beta,\delta, s)$ as given in Theorem~\ref{sing-prob2}. Furthermore,
\begin{equation*}
\Vert u\Vert_{C^\gamma(\overline{\Omega})}\leq C=C(\Vert h\Vert_{C_0(\overline{\Omega})},\beta,\delta,s,\lambda).
\end{equation*}
\end{lemm}
{\bf Proof of Lemma~\ref{lemmaA.1}.}\\
We observe as in the proof of Proposition~\ref{minimal_branch} that $\underline{u}_\lambda$ and $\underline{u}_\lambda+MU$, with $\underline{u}_\lambda$ and $U$ defined in \eqref{strictsub} and \eqref{supersol} and $M=\Vert h\Vert_{C_0(\overline{\Omega})}$, are subsolution and supersolution respectively to \eqref{eqA1}. Thus, using the comparison principle, we have that
\begin{equation*}
\underline{u}_\lambda\leq u\leq \underline{u}_\lambda+MU\mbox{ and }\Vert u\Vert_{C_{\phi_{\delta,\beta}(\Omega)}}\leq C(\Vert h\Vert_{C_0(\overline{\Omega})},\beta,\delta,s,\lambda).
\end{equation*}
Then, from $\frac{\lambda K(x)}{u^\delta}+h\leq C_1\phi_{\delta,\beta}^{\gamma-2s}$ for a constant $C_1>0$ and using Theorem~\ref{sing-prob2}, 
\begin{equation*}
\Vert u\Vert_{C^\gamma(\R^N)}\leq C(\Vert h\Vert_{C_0(\overline{\Omega})},\beta,\delta,s,\lambda).
\end{equation*}
\qed

From  Theorem\ref{sing-prob}, \cite[Proposition~1.2.9]{Abatangelo},\eqref{eq1.3} and the comparison principle, we establish the following lemma:
\begin{lemm}\label{lemmaA.2}
Let $a$ be a nonnegative and continuous function on $\Omega$. Let $u\in C_0(\overline{\Omega})$ be such that for $c\in\R$ and $\nu\in[0,2s)$:
\begin{equation*}
(-\Delta)^s u+a u\leq \frac{c}{d(x)^\nu}\mbox{ in }\Omega.
\end{equation*}
Then, for some constant $C>0$ independent of $u$
\begin{eqnarray*}
\left\{\begin{array}{lll}
& u\leq Ccd(x)^{s}\mbox{ if }\nu<s,\\
& u\leq Ccd(x)^s\ln\left(\frac{D}{d(x)}\right)\mbox{ with } D>{\rm diam}(\Omega)\mbox{ if }\nu=s,\\
& u\leq Ccd(x)^{2s-\nu}\mbox { if }\nu>s.
\end{array}\right.
\end{eqnarray*}
Similarly, if $u$ verifies
\begin{equation*}
(-\Delta)^s u+a u\leq \frac{c}{d(x)^s\ln^\alpha\left(\frac{D}{d(x)}\right)}\mbox{ in }\Omega \quad\mbox{with }0\leq \alpha<1,
\end{equation*}
then, $u\leq Ccd(x)^s\ln^{1-\alpha}\left(\frac{D}{d(x)}\right)$ for $C>0$ large enough and independent of $u$.
\end{lemm}
We now deal with the regularity of the function $A$. We first prove the continuity of $A$:
\begin{lemm}\label{lemmaA.3}
The map $A$ is continuous on $\R^+\times C_0(\overline{\Omega})$.
\end{lemm}
{\bf Proof of Lemma~\ref{lemmaA.3}.}\\
Let $h, h_\epsilon\in C_0(\overline{\Omega})$, $\lambda\in \R^+$, $\eta\in\R$, $A(\lambda,h)=u$ and $A(\lambda+\eta,h+h_\epsilon)=u_{\eta,\epsilon}$.
Then, we have
\begin{equation*}
(-\Delta)^s(u_{\eta,\epsilon})-\frac{(\lambda+\eta)K(x)}{u_{\eta,\epsilon}^\delta}=h+h_\epsilon\mbox{ in }\Omega\mbox{ and}
\end{equation*}
\begin{equation*}
(-\Delta)^s u-\lambda {K(x)}{u^\delta}=h\mbox{ in }\Omega.
\end{equation*}
Thus, for some $\theta\in(0,1)$
\begin{equation*}
(-\Delta)^s(u_{\eta,\epsilon}-u)+\frac{\lambda\delta K(x)}{(u+\theta(u_{\eta,\epsilon}-u))^{\delta+1}}(u_{\eta,\epsilon}-u)=h_\epsilon+\frac{\eta K(x)}{u_{\eta,\epsilon}^\delta}
\end{equation*}
From Lemma~\ref{lemmaA.2},  we get
\begin{equation*}
\vert u_{\eta,\epsilon}-u\vert \leq  O(\Vert h_\epsilon\Vert_{L^\infty(\Omega)}+\eta)\phi_{\delta,\beta}.
\end{equation*}
Therefore,
\begin{equation*}
\Vert u_{\eta,\epsilon}-u\Vert_{C_{\phi_{\delta,\beta}}(\Omega)}\to 0\mbox{ and  then } \Vert u_{\eta,\epsilon}-u\Vert_{C_0(\overline{\Omega})} \to 0\mbox{ as } (\Vert h_\epsilon\Vert_{L^\infty(\Omega)}+\vert\eta\vert)\to 0.
\end{equation*}
\qed

Finally we have:
\begin{lemm}\label{lemmaA.4}
$A$ is $C^2$ on $\R^+\times C_0(\overline{\Omega})$.
\end{lemm}
{\bf Proof of Lemma~\ref{lemmaA.4}.}\\
We first show that  for any $\lambda>0$, the map $C_0(\overline{\Omega})\ni u\mapsto A(\lambda,u)\in C_0(\overline{\Omega})$ is G\^ateaux differentiable. Let $\phi, h\in C_0(\overline{\Omega})$, $t\in\R$. We define $u_t=A(\lambda, h+t\phi)$ for $t>0$ and $u=A(\lambda, h)$. Then, we obtain
\begin{equation*}
(-\Delta)^s\left(\frac{u_t-u}{t}\right)+\frac{1}{t}\left(\frac{\lambda K(x)}{u^\delta}-\frac{\lambda K(x)}{u_t^\delta}\right)=\phi
\end{equation*}
which implies that for some $\theta\in(0,1)$
\begin{equation*}
(-\Delta)^s\left(\frac{u_t-u}{t}\right)+\frac{\lambda K(x)\delta}{(u+\theta(u_t-u))^{\delta+1}}\left(\frac{u_t-u}{t}\right)=\phi.
\end{equation*}
From Lemma~\ref{lemmaA.2}, $u+\theta(u_t-u)\geq c\phi_{\delta,\beta}$ with $c>0$ independent of $t$ and from Lemma~\ref{lemmaA.3} $\frac{u_t-u}{t}$ is bounded in $C_{\phi_{\delta,\beta}}(\Omega)$ and in $\tilde{H}^s(\Omega)$. Therefore, from Theorem~\ref{sing-prob2}, $\frac{u_t-u}{t}$ is bounded in $C^\gamma(\R^N)$, with $\gamma=\gamma(\delta,\beta, s)\in (0,1)$. Therefore,
\begin{equation*}
\frac{u_t-u}{t}\to v\mbox{ in }C_0(\overline{\Omega}) \mbox{ as } t\to 0^+
\end{equation*}
where $v$ satisfies
\begin{equation*}
(-\Delta)^s v+\frac{\lambda\delta K(x)}{u^{\delta+1}}v=\phi\mbox{ in }\Omega\mbox{ and where } u=A(\lambda, h).
\end{equation*}
We can also show easily that the map $\phi\mapsto v$ is continuous in $C_0(\overline{\Omega})$. This proves the G\^ateaux differentiability of $u\to A(\lambda,u)$ and $\partial_2A(\lambda,h)(\phi)=v$. Next, we prove the Frechet-differentiability of $u\to A(\lambda,u)$ (with $\lambda\in\R^+$ fixed). For $\phi\in C_0(\overline{\Omega})$, we define $u_\phi\eqdef A(\lambda, h+\phi)$ and recall $v=\partial_2A(\lambda,h)(\phi)$. Then, we get for suitable $\theta_0,\theta_1\in (0,1)$,
\begin{eqnarray*}
&(-\Delta)^s\left(\frac{u_\phi-u-v}{\Vert \phi\Vert_{L^\infty(\Omega)}}\right)+\frac{\lambda\delta K(x)}{(u+\theta_0(u_\phi-u))^{\delta+1}}\left(\frac{u_\phi-u-v}{\Vert \phi\Vert_{L^\infty(\Omega)}}\right)=\left(\frac{\lambda\delta K(x)}{u^{\delta+1}}-\frac{\lambda\delta K(x)}{(u+\theta_0(u_\phi-u))^{\delta+1}}\right)\frac{v}{\Vert \phi\Vert_{L^\infty(\Omega)}}\\
&=\frac{\lambda\delta(\delta+1)\theta_0 K(x)(u_\phi-u)v}{(u+\theta_1(u_\phi-u))^{\delta+2}\Vert \phi\Vert_{L^\infty(\Omega)}}\mbox{ in }\Omega.
\end{eqnarray*}
For some $\theta_2\in (0,1)$, we have also that
\begin{equation*}
(-\Delta)^s(u_\phi-u)+\frac{\lambda\delta K(x)(u_\phi-u)}{(u+\theta_2(u_\phi-u))^{\delta +1}}=\phi\mbox{ in }\Omega
\end{equation*}
and from Lemma~\ref{lemmaA.2}
\begin{equation*}
\vert u_\phi-u\vert\leq C\Vert \phi\Vert_{L^\infty(\Omega)}\phi_{\delta,\beta}.
\end{equation*}
Then, 
\begin{equation*}
(-\Delta)^s\left(\frac{u_\phi-u-v}{\Vert\phi\Vert_{L^\infty(\Omega)}}\right)+\frac{\lambda\delta K(x)}{(u+\theta_0(u_\phi-u))^{\delta+1}}\left(\frac{u_\phi-u-v}{\Vert\phi\Vert_{L^\infty(\Omega)}}\right)=\frac{o_\phi(1)}{u^\delta}.
\end{equation*}
Therefore,
\begin{equation*}
\displaystyle\frac{\Vert u_\phi-u-v\Vert_{L^\infty(\Omega)}}{\Vert\phi\Vert_{L^\infty(\Omega)}}\to 0\mbox{ as }\Vert\phi\Vert_{L^\infty(\Omega)}\to 0.
\end{equation*}
This proves the Frechet differentiablity of $u\mapsto A(\lambda,u)$. We now prove that this function is $C^1$ on $C_0(\overline{\Omega})$. We need to prove that $h\to \partial_2 A(\lambda, h)$ is continuous. Let $(h_n)_{n\in\N}\subset C_0(\overline{\Omega})$ such that $h_n\to h$ in $C_0(\overline{\Omega})$. Then,
\begin{eqnarray*}
\Vert\partial_2 A(\lambda, h_n)-\partial_2 A(\lambda, h)\Vert=\displaystyle\sup_{0\not\equiv\phi\in C_0(\overline{\Omega})}\frac{\Vert(\partial_2 A(\lambda, h_n)-\partial_2 A(\lambda, h))(\phi)\Vert_{L^\infty(\Omega)}}{\Vert\phi\Vert_{L^\infty(\Omega)}}.
\end{eqnarray*}
Setting 
\begin{equation*}
v\eqdef\partial_2 A(\lambda, h))(\phi)\mbox{ and }v_n\eqdef\partial_2 A(\lambda, h_n)(\phi),
\end{equation*}
we have
\begin{eqnarray*}
&(-\Delta)^s\left(\frac{v-v_n}{\Vert\phi\Vert_{L^\infty(\Omega)}}\right)+\frac{\lambda\delta K(x)}{A(\lambda,h)^{\delta+1}}\left(\frac{v-v_n}{\Vert\phi\Vert_{L^\infty(\Omega)}}\right)=\left(\frac{\lambda\delta K(x)}{A(\lambda,h)^{\delta+1}}-\frac{\lambda\delta K(x)}{A(\lambda,h_n)^{\delta+1}}\right)\frac{v_n}{\Vert\phi\Vert_{L^\infty(\Omega)}}\\
&=\frac{o_{\Vert h_n-h\Vert_{L^\infty(\Omega)}}(1)}{A(\lambda,h)^{\delta+1}}.
\end{eqnarray*}
Therefore using Lemma~\ref{lemmaA.2}, we obtain that
\begin{equation*}
\frac{\Vert v-v_n\Vert_{L^\infty(\Omega)}}{\Vert\phi\Vert_{L^\infty(\Omega)}}\to 0\mbox{ as }\Vert h_n-h\Vert_{L^\infty(\Omega)}\to 0.
\end{equation*}
This completes the proof of  the $C^1$ regularity of $u\to A(\lambda, u)$. Similarly, we can prove that $(\lambda, u)\to A(\lambda,u)$ is $C^2$ withe following continuous partial derivatives: $\partial_2 A(\lambda, h)(\phi)=v$, $\partial_1 A(\lambda, h)=w_1$, $\partial^2_{11}A(\lambda, h)=w_{11}$, $\partial^2_{22}A(\lambda, h)(\phi,\psi)=w_{22}$, $\partial^2_{12}A(\lambda, h)(1,\phi)=w_{12}\in C_{\phi_{\delta,\beta}}(\Omega)\cap \tilde{H}^s(\Omega)$ satisfy
\begin{equation*}
(-\Delta)^sw_1+\frac{\lambda\delta K(x)}{A(\lambda,h)^{\delta+1}}w_1=\frac{K(x)}{A(\lambda,h)^\delta}\mbox{ in }\Omega,
\end{equation*}
\begin{equation*}
(-\Delta )^sw_{11}+\frac{\lambda\delta K(x)}{A(\lambda,h)^{\delta+1}}w_{11}=\frac{\delta(\delta+1)K(x)w_1^2}{A(\lambda,h)^{\delta+2}}-\frac{2K(x)\delta w_1}{A(\lambda,h)^{\delta+1}}\mbox{ in }\Omega,
\end{equation*}
\begin{equation*}
(-\Delta)^s w_{12}+\frac{\lambda\delta K(x)}{A(\lambda,h)^{\delta+1}}w_{12}=\frac{\delta(\delta+1)K(x)w_1v}{A(\lambda,h)^{\delta+2}}-\frac{\delta K(x)v}{A(\lambda,h)^{\delta+1}}\mbox{ in }\Omega,
\end{equation*}
\begin{equation*}
(-\Delta)^s w_{22}+\frac{\lambda\delta K(x)}{A(\lambda,h)^{\delta+1}}w_{22}=\frac{\delta(\delta+1)K(x)v^2}{A(\lambda,h)^{\delta+2}}\mbox{ in }\Omega
\end{equation*}
for any $\lambda>0$, $h,\phi,\psi\in C_0(\overline{\Omega})$. We omit the details here.
\qed

{\bf Acknowledgement:} Adimurthi and J. Giacomoni were funded by IFCAM (Indo-French Centre for Applied Mathematics, UMI CNRS 3494) under the project ``Singular phenomena in reaction diffusion equations and in conservation laws". J. Giacomoni thanks TIFR CAM Bangalore for its kind hospitality and the nice scientific atmosphere during his visit in fall 2016. S. Santra acknowledges support from LMAP UMR UPPA-CNRS 5142, Universit\'e de Pau et des Pays de l'Adour.


\begin{thebibliography}{99}

\bibitem{Abatangelo} N. Abatangelo, {\it Large S-harmonic functions and boundary blow-up solutions for the fractional Laplacian}, Discrete Contin. Dyn. Syst., 35(12) (2015), 5555--5607.

\bibitem{AJ} Adimurthi and J. Giacomoni, {\it  Multiplicity of positive solutions for a singular and critical elliptic
problem in $\mathbb R^2$}, Commun. Contemp. Math., 8(5) (2006), 621--656.

\bibitem{da} D. Applebaum, {\it L$\acute{e}$vy processes-from probability to finance and quantum groups}, Notices Amer. Math. Soc.,
51 (2004), 1336--1347.

\bibitem{BaGi} K. Bal and J. Giacomoni, {\it A remark on symmetry of solutions to singular equations and applications}, Eleventh International Conference Zaragoza-Pau on Applied Mathematics and Statistics, Monogr. Mat. Garc\'ia Galdeano, 37 (2012), 25--35.

\bibitem{peral-al}  B. Barrios, I. De Bonis,  M. Medina and  I. Peral, {\it  Semilinear problems for the fractional laplacian with a singular nonlinearity}, Open Math., 13 (2015), 390--407.


\bibitem{BoGiPr} B. Bougherara, J. Giacomoni and S. Prashanth, {\it Analytic global bifurcation and infinite turning points for very singular problems}, Calc. Var. Partial Differential Equations, 52(3-4) (2015), 829--856.

\bibitem{BrMe} H. Brezis and F. Merle, {\it Uniform estimates and blow-up behavior
  for solutions of $-\Delta u=V(x)e^u$ in two dimensions}, Comm. P.D.E., 16 (1991), 1223--1253.

\bibitem{Bucur} C. Bucur, {\it Some observations on the Green function for the ball in the fractional Laplace framework}, Commun. Pure Appl. Anal., 15(2) (2016), 657--699.

\bibitem{Buffoni-Toland-book} B. Buffoni and J.~F. Toland,  {\sl Analytic theory of global bifurcation. An introduction},
 Princeton Series in Applied Mathematics, Princeton University Press, Princeton, NJ, 2003.

\bibitem{CS} L. Caffarelli, L. Silvestre, {\it An extension problem related to the fractional Laplacian}, Comm.
    Partial Differential Equations, 32 (2007), 1245--1260.

\bibitem{Chen-Li-Arma} W. Chen and C. Li, {\it A priori estimates for solutions to nonlinear elliptic equations}, Arch. Rational Mech. Anal., 122(2) (1993), 145--157.

\bibitem{Chen-Li-Arxiv} W. Chen, C. Li and Y. Li, {\it A direct method of moving planes for the fractional Laplacian}, Advances in Math. 308, (2017),  404--437.

\bibitem{Chen-Song} Zhen-Qing Chen and Renming Song, {\it Estimates on Green functions and Poissons kernels for symmetric stable processes}, Math. Ann., 312 (1998), 465--501.

\bibitem{ChSoKi} Zhen-Qing Chen, Panki Kim and Renming Song, {\it  Heat kernel estimates for the Dirichlet fractional Laplacian}, J. Eur. Math. Soc., 12 (2010), 1307--1329.

\bibitem{Cr-Ra} M.~G.  Crandall and P.~H. Rabinowitz, {\it Bifurcation from simple eigenvalues}, J. Funct. Anal., 8 (1971), 321--340.

\bibitem{Crandall-Rabinowitz-arma} M.~G. Crandall and P.~H. Rabinowitz, {\it Bifurcation Perturbation of Simple Eigenvalues and linearized Stability}, Arch. Ration. Mech. Anal., 52 (1973), 161--180.

\bibitem{Cr-Ra-Ta} M. G. Crandall,  P. H. Rabinowitz and L. Tartar, {\it On a Dirichlet problem with a singular nonlinearity}, Comm.  Partial Differential Equations,  2 (1977), 193--222.

\bibitem{Da Lio-Martinazzi-Riviere} F. Da Lio, L. Martinazzi and T. Rivi\`ere, {\it Blow-up analysis of a nonlocal Liouville-type equation}, Anal. PDE, 8(7), (2015), 1757--1805.

\bibitem{Dancer} E.~N. Dancer, {\it Bifurcation theory for analytic operators}, Proc. Lond. Math. Soc., XXVI (1973), 359--384.

\bibitem{Da2} E.~N. Dancer, {\it Real analyticity and non degeneracy}, Math. Ann., 325(2) (2003), 369--392.

\bibitem{DhGiPrSa} R. Dhanya, J. Giacomoni,  S. Prashanth, and K. Saoudi, {\it Global bifurcation and local multiplicity results for elliptic equations with singular nonlinearity of super exponential growth in $\mathbb{R}^2$},  Adv.  Differential Equations, 17(3-4) (2012), 369-400.

\bibitem{DiMoOs} J.~I. D\'{\i}az, J.~M. Morel and L. Oswald, {\it An elliptic equation with singular
nonlinearity}, Commun. Partial Differential Equations, 12 (1987), 1333--1344.

\bibitem{Franzina-Palatucci} G. Franzina and G. Palatucci, {\it Fractional $p$-eigenvalues}, Riv. math. Univ. Parma (N.S), 5(2), 373--386.

\bibitem{FuMa} W. Fulks and J.~S. Maybee, {\it A singular nonlinear equation}, Osaka J. Math., 12, (1960), 1--19.

\bibitem{GaJu} I.~M. Gamba and A. Jungel, {\it Positive solutions to a singular second and third
order differential equations for quantum fluids}, Arch. Ration. Mech. Anal., 156 (2001), 183--203.

\bibitem{GhRa} M. Ghergu and V. R\u{a}dulescu, {\it Singular elliptic problems: bifurcation and asymptotic analysis}, Oxford University Press, 2008.

\bibitem{GhRa2} M. Ghergu and  V. R\u{a}dulescu, {\it Multiparameter bifurcation and asymptotics for the singular Lane-Emden-Fowler equation with a convection term}, Proceedings of the Royal Society of Edinburgh: Section A (Mathematics), 135 (2005), 61--84.

\bibitem{GiPaSe} J. Giacomoni, P.~K. Mishra and K. Sreenadh, {\it Fractional elliptic equations with critical exponential nonlinearities},  Advances in Nonlinear Analysis, 5(1) (2016), 57--74.

\bibitem{TuGiSe} J. Giacomoni, T. Mukherjee and K. Sreenadh, {\it Positive solutions of fractional elliptic equation with critical and singular nonlinearty}, to appear in Adv. Nonlinear Anal., DOI: 10.1515/anona-2016-0113.

\bibitem{GiPrWa} J. Giacomoni, S.K. Prashanth and G. Warnault, {\it Existence and global analytic bifurcation for singular biharmonic equation with Navier boundary condition}, Proc. Amer. Math. Soc., 145(1) (2017), 151--164.

\bibitem{Gidas-Spruck} B. Gidas et J. Spruck, {\it A priori bounds for positive solutions of nonlinear equations}, Comm. P.D.E., 6 (1981), 883--901.

\bibitem{Gr} P. Grisvard, {\sl Elliptic Problems in nonsmooth domains}, Monogr. Stud. Math., vol. 24, Pitman (Advances Publishing Program), Boston, MA, 1985.
\bibitem{HeMa} J. Hern\'andez and F.~J. Mancebo, {\sl Singular elliptic and parabolic equations}, Handbook of Differential Equations,   3 (2006), 317--400.

\bibitem{Figueiredo-Lions_Nussbaum-JMPA} D.~G. de Figueiredo, P.~L. Lions and R.~D. Nussbaum, {\it Apriori estimates and existence of positive solutions of semilinear elliptic equations}, J. math. Pures Appl. (9), 61(1) (1982), 41--63.

\bibitem{Martinazzi} L. Martinazzi, {\it Fractional Adams-Moser-Trudinger type inequalities}, Nonlinear Anal., 127 (2015), 263--275.

\bibitem{Ra-JFA} P.~H. Rabinowitz, {\it Some global results for nonlinear eigenvalue problems}, J. Functional Analysis, 7 (1971), 487--513.

\bibitem{Ros-oton-serra-JMPA} X. Ros-Oton and J. Serra, {\it The Dirichlet problem for the fractional Laplacian: Regularity up to the boundary}, J. Math. Pures Appl., 101 (2014), 275--302.

\bibitem{Ros-othon-serra-CVPDE} X. Ros-Oton and J. Serra, {\it  The extremal solution for the fractional Laplacian}, Calc. Var. Partial Differential Equations, 50(3-4) (2014), 723--750.

\bibitem{Silvestre-CPAM} L. Silvestre, {\it Regularity of the obstacle problem for a fractional power of the Laplace operator}, Comm. Pure Appl. Math., 60(1 (2007), 67--112.

\bibitem{Tr} H. Triebel, {\sl Interpolation Theory, Function Spaces, Differential Operators}, second ed., Johann Ambrosius Barth, Heidelberg, 1995.
\end{thebibliography}
\end{document}